\documentclass[a4paper,11pt]{amsart}

\usepackage{a4wide}
\usepackage{amsmath,amssymb,amsthm}
\usepackage{amsfonts,amsthm,amsmath,amssymb,amscd,color}
\usepackage{enumerate}
\usepackage{parskip}
\usepackage{graphicx}
\usepackage{babel}

\newcommand{\RR}{\mathbb{R}}

\newcommand{\NN}{\mathbb{N}}

\newcommand{\eps}{\varepsilon}

\newcommand{\CE}{{\mathcal{E}}}

\newcommand{\CM}{{\mathcal{M}}}

\newcommand{\BL}{\mathrm{BL}}

\newcommand{\FM}{\mathrm{FM}}

\newcommand{\supp}{\mathrm{supp}}
\newcommand{\Lip}{\mathrm{Lip}}

\newtheorem{thrm}{Theorem}[section]

\newtheorem{prop}{Proposition}[section]
\newtheorem{lemma}{Lemma}[section]
\newtheorem{clry}{Corollary}[section]

\newtheorem{example}{Example}[section]
\newtheorem{remark}{Remark}[section]

\newcommand{\smfrac}[2]{\mbox{$\frac{#1}{#2}$}}
\newcommand{\pair}[2]{\left\langle #1 , #2 \right\rangle}
\newcommand{\mathbbm}[1]{{#1\!\!#1}}
\newcommand{\ind}{{\mathbbm{1}}}
\newcommand{\Mol}{{\mathcal{M}ol}}
\newcommand{\ext}{{\mathrm{ext}}}

\begin{document}
	\title[Norming and dense sets of extreme points]{Norming and dense sets of extreme points of the unit ball in spaces of bounded Lipschitz functions}
	
	
	\author{Sander C. Hille}
	\address{Mathematical Institute, Leiden University, P.O. Box 9512, 2300 RA Leiden, The Netherlands}
	\email{shille@math.leidenuniv.nl}
	
	\author{Esm\'ee S. Theewis}
	\address{Delft Institute of Applied Mathematics, Delft University of Technology, P.O. Box 5031, 2600 GA Delft, The Netherlands, (ET)}
	\email{e.s.theewis@tudelft.nl}
	
	\date\today
	
	
	\begin{abstract}
		On spaces of finite signed Borel measures on a metric space one has introduced the Fortet-Mourier and Dudley norms, by embedding the measures into the dual space of the Banach space of bounded Lipschitz functions, equipped with different -- but equivalent -- norms: the FM-norm and the BL-norm, respectively. The norm of such a measure is then obtained by maximising the value of the measure when applied by integration to extremal functions of the unit ball. We introduce suitable subsets of extremal functions in both cases, essentially by a construction by means of McShane's Lipschitz extension operator for Lipschitz functions on finite metric spaces, that are shown to be norming. As a consequence, we obtain that these sets are weak-star dense in the full set of extreme points of the unit ball. If the underlying metric space is connected, then we have found a larger set of extremal functions for the BL-norm, similar to such a set that was defined previously by J. Johnson for the FM-norm when the underlying space was compact. This set is then also weak-star dense in the extremal functions. This may open an avenue to obtaining computational approaches for the Dudley norm on signed Borel measures.
	\end{abstract}
	
	\maketitle

\section{Introduction}

Consider a metric space $(S,d)$ and the vector space $\Lip(S,d)$ of Lipschitz functions for the metric $d$. We shall fix the metric throughout this study, so it will be omitted in notation. The Lipschitz constant of $f\in\Lip(S)$ is denoted by $|f|_L$. Banach spaces of Lipschitz functions are of interest, for example because they induce norms on finite signed measures on $S$, by viewing measures as functionals on theses Banach spaces by integration (see e.g. \cite{Hille-Worm:2009,Gwiazda-Thieme_ea:2018,Lasota-Myjak-Szarek:2002}). The norm is then the restriction of the dual norm to the embedded space of measures.

Here, we are interested in the space of bounded Lipschitz functions $\BL(S)$, with the norm
\begin{equation}\label{eq:def BL-norm}
\|f\|_\BL := \|f\|_\infty + |f|_L.
\end{equation}
It is a Banach space. Let $B_\BL^S$ denote its closed unit ball and let $\CM(S)$ denote the finite signed Borel measures on $S$. Embedding of a finite signed Borel measure $\mu$ on a complete separable metric space $(S,d)$ results in the {\it Dudley norm} or {\it dual bounded Lipschitz norm}
\begin{equation}\label{eq:dual BL-norm}
	\|\mu\|_\BL^* := \sup\bigl\{ \int_S f\, d\mu\colon f\in B_\BL^S \bigr\}
\end{equation}
on $\CM(S)$. It metrizes the weak topology on the finite  positive measures $\CM^+(S)$, which is defined by pairing with bounded continuous functions by integration (cf. \cite{Dudley:1966,Dudley:1974,Hille-Worm:2009} and also \cite{Gwiazda-Thieme_ea:2018}) and is commonly used in probability theory.

In mathematical theory, estimates of the norm are often sufficient. In various applications though, e.g. numerical analysis one wants to be able to compute arbitrarily good approximations of the norm $\|\mu\|_\BL^*$, for suitable classes of measure $\mu$, to be able to estimate error. To that end, one can show that the supremum in \eqref{eq:dual BL-norm} is actually attained at an extreme point of $B_\BL^S$ (see Proposition \ref{prop:norm attained in extreme point}). Therefore, it is of interest to determine a set of extreme points of $B^S_\BL$ to which the suppremum in \eqref{eq:dual BL-norm} can be restricted for determining the $\|\cdot\|_\BL^*$ norm of measures $\mu\in\CM(S)$. We call such a set `{\it norming}'.

In \cite{Hille-Theewis:2022}, we succeeded in achieving the objective of deriving particular explicit expressions and computational approaches to the so-called {\it Fortet-Mourier norm} on $\CM(S)$ (see e.g. \cite{Lasota-Myjak-Szarek:2002}),
\begin{equation}\label{eq:FM norm}
\|\mu\|_\FM^* := \sup \{ \int_S f\, d\mu\colon f\in\BL(S),\ \|f\|_\infty\leq 1,\ |f|_L\leq 1\},
\end{equation}
for any Polish space $S$. This norm is associated to the norm 
\[
\|f\|_\FM\ :=\ \max\bigl(\|f\|_\infty,|f|_L\bigr).
\]
on $\BL(S)$, which is equivalent to $\|\cdot\|_\BL$. Hence, $\|\cdot\|_\FM^*$ is equivalent to $\|\cdot\|_\BL^*$. Previously, a computational approach for determining $\|\mu\|_\FM$ was available only for very particular cases, e.g. $S=[0,1]$ and $\mu$ finitely supported, see \cite{Jablonski-Marciniak-Czochra:2013}.

We could do so, because a characterisation of the extreme points of the unit ball $B_\FM^S$ of $(\BL(S),\|\cdot\|_\FM)$ had been provided by Farmer \cite{Farmer:1994} in general (though for the seemingly different metric dual space $(\Lip_0(S),|\cdot|_L)$ introduced by Lindenstrauss \cite{Lindenstrauss:1964}, that turned out to be isometrically isomorphic to $(\BL(S),\|\cdot\|_\FM)$, see Appendix \ref{app:isomormism FM-case}), while a few particular cases had been covered before by Roy \cite{Roy:1968}, Rao \& Roy \cite{Rao-Roy:1970}, and Rolewicz \cite{Rolewicz:1986}, in particular $S=[0,1]$. See also \cite{Cobzas:1989,Smarewski:1997}. Moreover, Johnson \cite{Johnson:1975} had described a dense subset of the extreme points of $B^S_\FM$ for compact and connected $S$ that suffices to compute the norm \eqref{eq:FM norm}. In particular the latter, provided the inspiration for the results obtained in \cite{Hille-Theewis:2022}. The techniques of proof that we employed do not work for $\|\cdot\|_\BL$ though, essentially because one cannot change the $\|\cdot\|_\infty$-norm of a function in the boundary of $B^S_\BL$ while staying within $B^S_\BL$, without changing also its Lipschitz constant. For $B^S_\FM$ these can be changed independently.

To our knowledge, no results existed in the literature regarding the characterisation of the extreme points of $B^S_\BL$ or the description of a dense subset of such points, similar to Johnson's, suitable for computing the norm $\|\cdot\|^*_\BL$. Here we provide two dense and norming sets: a `small' set $E^S_\BL$, constructed through Lipschitz extension operators and a `Johnson-like' set $J_\BL^S$ that consists of extreme points when $S$ is connected. 
Our main results are Theorem \ref{thrm:main result dense set extreme points} and its corollaries, Corollary \ref{clry:E is norming} and Corollary \ref{clry:density Johnson-like set BL}, all for general metrics spaces $(S,d)$. No connectedness, nor compactness condition needs to be imposed on $S$ for the `small' set $E^S_\BL$ to consist of extreme points. The set $J^S_\FM$ and the new set $J^S_\BL$ do consist of extreme points if $S$ is connected. No compactness of $S$ is required (cf. Proposition \ref{E_BL is subset}).

The description of $E^S_\BL$ and the proof of the theorem requires revisiting the classical McShane extension operator (cf. \cite{McShane:1934}, Theorem 1) that extends Lipschitz functions defined on a subset $P$ of $S$ (with induced metric) to the whole space $S$ without increasing its Lipschitz constant. Upon further inspection, this non-linear operator, denoted by $\CE_P^{S,0}$, turns out to preserve the Lipschitz constant. We give an overview of its properties and that of a slightly modified Lipschitz extension operator $\CE_P^S$ that plays a pivotal role in the Metric Tietze Extension Theorem (Theorem \ref{thrm:metric Tietze extension}), which is central in the argumentation in this paper. Because of this centrality and because a full proof of this results is not readily available in the literature (although reference and a sketch of proof is found in \cite{Weaver:1999}, Theorem 1.5.6), we provide a full proof in Appendix \ref{App:proof metric Tietze extension thrm}.

A particular result that is worth noting in this context, is Theorem \ref{thrm:extension non-trivial extreme points}, that states that the extension operator $\CE_P^S$ maps non-trivial extreme points of $B^P_\BL$ into non-trivial extreme points of $B^S_\BL$ (non-triviality will be specified in Section \ref{sec:basic properties extreme points}). This was known for $(\Lip_0(S),|\cdot|_L)$ -- hence, in disguise for $(\BL(S),
\|\cdot\|_\FM)$ -- by Farmer \cite{Farmer:1994} (see Proposition 3). His proof relied on his characterisation of extreme points. Our proof for $(\BL(S), \|\cdot\|_\BL)$ can be modified slightly, then covering also the case of $\|\cdot\|_\FM$, without the need of a characterisation of the extreme points. This proof is sketched in Appendix \ref{App:Norming sets FM norm}.

Moreover, on our way towards the main result, we established various auxiliary results on extreme points of $B^S_\BL$ and $B^S_\FM$. We found several of these of sufficiently general interest to collect them in Section \ref{sec:basic properties extreme points}.

\section{Notation and preliminary results}

Throughout, $(S,d)$ will be a metric space. The main object of interest in this paper are the unit balls in $\BL(S)$ and their convex structure, mainly for the norm $\|\cdot\|_\BL$, but occasionally also for $\|\cdot\|_\FM$, which will be denoted by  
\[
B^S_\bullet := \bigl\{ f\in \BL(S): \|f\|_\bullet\leq 1\bigr\}\qquad\mbox{with}\ \bullet=\BL,\FM.
\]

A measure $\mu\in\CM(S)$ is called (strictly) {\it separable} if there exists a separable set $S_0\subset S$, such that $\mu$ is concentrated on $S_0$, i.e. $|\mu|(S_0)=|\mu|(S)$, where $|\mu|=\mu^++\mu^-$ is the total variation measure associated to $\mu$. We denote the linear subspace of separable measures in $\CM(S)$ by $\CM_s(S)$. If $S$ is separable, then $\CM(S)=\CM_s(S)$. The map $\mu\mapsto I_\mu:\CM(S)\to\BL(S)^*$, defined by
\begin{equation}\label{eq:pairing functions-measures}
	I_\mu(f) := \langle\mu,f\rangle := \int_S f\,d\mu,\qquad f\in\BL(S),
\end{equation}
is injective, because on a metric space, any $\mu\in\CM(S)$ is regular (cf. \cite{Bogachev-II:2007}, Theorem 7.7.1, p.70). This embedding makes it natural to equip $\CM(S)$ with the restriction of the dual norm $\|\cdot\|^*_\FM$ or $\|\cdot\|^*_\BL$ on $\BL(S)^*$. This is the idea behind the introduction of the Dudley norm \eqref{eq:dual BL-norm} or the Fortet-Mourier norm \eqref{eq:FM norm} on measures. Thus, 
\begin{equation}\label{eq:expression dual norm}
	\|\mu\|_\bullet^* := \sup_{f\in B^S_\bullet} \bigl|\langle\mu,f\rangle\bigr| = \sup_{f\in B^S_\bullet} \langle\mu,f\rangle,\qquad \bullet=\FM,\BL.
\end{equation}
These two norms on measures are equivalent. $\CM(S)_\BL$ denotes the space $\CM(S)$ equipped with the $\|\cdot\|_\bullet^*$-norm topology. $\overline{\CM}(S)_\BL$ and $\overline{\CM}_s(S)_\BL$ are the completions of $\CM(S)$ and $\CM_s(S)$, viewed as closure in $(\BL(S)^*,\|\cdot\|_\BL^*)$.

Following Pachl \cite{Pachl:2013}, any finite real linear combination of Dirac measures is called a {\it molecular measure}:
\begin{equation}
	\Mol(S) := \mathrm{span}_\RR\bigl\{ \delta_x: x\in S\bigr\}.
\end{equation}
Clearly, $\Mol(S)\subset \CM_s(S)$. It is dense for $\|\cdot\|_\bullet^*$, ($\bullet=\FM, \BL$).

\begin{prop}\label{prop:dual Ms is BL}
	The map $\phi\mapsto f_\phi:\bigl(\Mol(S),\|\cdot\|_\bullet^*\bigr)^*\to \bigl(\BL(S),\|\cdot\|_\bullet^*\bigr)$, with $f_\phi(x):= \phi(\delta_x)$, $x\in S$, is a linear isometric isomorphism for $\bullet = \FM, \BL$.
\end{prop}

$\BL(S)$ is a Banach algebra for the norm $\|\cdot\|_\BL$ (\cite{Dudley:1966}, Lemma 3). Moreover, it is a vector lattice (or Riesz space) for the usual point-wise partial ordering `$\leq$' of functions. The lattice operations are given pointwise:
\begin{align*}
	\sup(f,g)(x)\ &=\ (f\vee g)(x)\ =\ \max\bigl( f(x),g(x)\bigr),\\
	\inf(f,g)(x)\ &=\ (f\wedge g)(x)\ =\ \min\bigl( f(x),g(x)\bigr),
\end{align*}
for any $f,g\in\BL(S)$. In particular, for $f\in\BL(S)$, $f^+=f\vee 0$ and $f^-= -(f\wedge 0)$ are in $\BL(S)$. One has
\[
|f\vee g|_L,\ |f\wedge g|_L\ \leq\ \max\bigl( |f|_L,|g|_L\bigr)
\]
(cf. \cite{Dudley:1966}), Lemma 4). $\|\cdot\|_\BL$ is not a Riesz norm.

Let $C$ be a convex set in a vector space $V$. $e\in C$ is an extreme point if there is no open line segment that contains $e$ and lies entirely in $C$. Equivalently, (see \cite{Day:1973}, p.101), $e\in C$ is an extreme point if and only if:
\begin{equation}
	e = \smfrac{1}{2}(x+y),\ x,y\in C\qquad \mbox{if and only if}\qquad x=y=e.
\end{equation}
The set of extreme points of $C$ is denoted by $\ext(C)$.

Proposition \ref{prop:dual Ms is BL} has the following consequence:
\begin{prop}\label{prop:norm attained in extreme point}
	Let $(S,d)$ be a metric space and $\mu\in\CM_s(S)$. Let $\bullet=\FM$ or $\BL$. Then there exists $f^\bullet_\mu\in\ext(B^S_\bullet)$ such that
	$\|\mu\|_\bullet^* = \langle\mu,f^\bullet_\mu\rangle$.
\end{prop}
\begin{proof}
	According to Proposition \ref{prop:dual Ms is BL} and the Banach-Alaoglu Theorem (cf. \cite{Conway:1990}, Theorem V.3.1, p.130), $B^S_\bullet$ is compact for the $\sigma(\BL(S),\overline{\Mol}(S)_\BL)$-topology. \cite{Bourbaki:topol}, Proposition 1, p. II.54 yields that the convex function $f\mapsto\langle\mu,f\rangle$ attains a maximum value at some $f^\bullet_\mu\in\ext(B^S_\bullet)$ of the weak*-compact convex set $B^S_\bullet$. Expression \eqref{eq:expression dual norm} completes the argument.
\end{proof}

The following alternative characterisation of extreme points will be the `workhorse' in our argumentation in proofs later on:
\begin{lemma}\label{lem:char extreme point}
	Let $C$ be a convex set in a vector space $V$. $e\in\ext(C)$ if and only if
	\begin{equation}\label{eq:cond extreme point}
		\bigl\{ x\in V:\ e+x\in C\ \mathrm{and}\ e-x\in C\bigr\} = \{0\}.
	\end{equation}
\end{lemma}
\begin{proof}
	If $e\in C$ is not an extreme point, then $e=\frac{1}{2}(x+y)$ for some $x,y\in C$, $x\neq y$. Put $z=\frac{1}{2}(x-y)$. Then $z\neq 0$, while $e+z = x\in C$ and $e-z=y\in C$. On the other hand, if \eqref{eq:cond extreme point} does not hold, then there exists $z\neq 0$ in $V$ such that $x:=e+z$ and $y:=e-z$ are in $C$. Then $x\neq y$, while $e=\smfrac{1}{2}(x+y)$. Therefore, $e$ cannot be an extreme point of $C$.
\end{proof}
Two immediate consequences of this lemma can be observed.

\begin{clry}\label{clry:extremes on sphere}
	Let $(X,\|\cdot\|)$ be a normed space, $R>0$ and put $B_R:= \{ x\in X: \|x\|\leq R\}$. Then $\ext(B_R)\subset\{x\in X: \|x\|=R\}$.
\end{clry}
\begin{proof}
	If $e\in \ext(B_R)$ and $\|e\|<R$, take $x\in X$ with $0<\|x\|\leq R-\|e\|$. Then $x\neq 0$, while $e+x$ and $e-x$ are in $B_R$. This contradicts Lemma \ref{lem:char extreme point}.
\end{proof}
A set $C$ in a vector space $V$ is called symmetric around 0, if $x\in C$ implies $-x\in C$. That is, $C=-C$.
\begin{clry}\label{clry:symmetric sets and extreme points}
	Let $C$ be a convex set in a vector space $V$. If $C$ is symmetric around 0, then $\ext(C)$ is symmetric around 0.
\end{clry}

We will need the following two functional analytic results, of which the second is called the K$^2$-M$^3$-R Theorem in \cite{Day:1973} for the contributions of Krein, Klee, Milman, Mazur, Minkowski and Rutman to various parts of the overall result, see \cite{Day:1973}, p.104 and also \cite{Klee:1957}, Theorem 1.1. It summarizes various characterisations on a subset of the extreme points of a compact convex set, such that this subset still `generates' this compact set as its closed convex hull.

\begin{lemma}\label{lem:weakstar dense subsets unit ball}
	Let $X$ be a Banach space and $F$ a dense subset of $X$. Let $B^*$ be the closed unit ball in $X^*$ and $A\subset B^*$. If $\|x\|_X = \sup\bigl\{ \phi(x)\colon \phi\in A\bigr\}$ for all $x\in F$, then the weak$^*$-closed convex hull of $A$ is $B^*$.
\end{lemma}
\begin{proof}
	See Lemma 1.1 in \cite{Johnson:1975}.
\end{proof}

\begin{thrm} (K$^2$-M$^3$-R Theorem)\label{thrm:K2-M3-R}
	Let $K$ be a convex compact set in a locally convex Hausdorff topological vector space $X$. Put $E=\ext(K)$ and $A\subset K$. Then the following are equivalent:
	\begin{enumerate}
		\item[({\it i})] The closed convex hull of $A$ is $K$;
		\item[({\it ii})] The closure of $A$ contains $E$;
		\item[({\it iii})] The closure of $A$ contains at least one point of each minimal facet of $K$;
		\item[({\it iv})] For each continuous linear functional $\phi$ on $X$, $\sup \phi(A) = \sup \phi(K)$.
	\end{enumerate}
\end{thrm}

\begin{proof}
	See \cite{Day:1973}, Section V.1.
\end{proof}

\section{Basic properties of extreme points of balls of Lipschitz functions}
\label{sec:basic properties extreme points}

Let $(S,d)$ be a metric space, without further assumptions (like completeness or separability). The geometry of the unit ball in $\BL(S)$ depends on the chosen norm. Hence, so does the set of extreme points. The $\FM$-norm and $\BL$-norm, though equivalent, introduce quite a different geometry. The argumentation required to exhibit extreme points in their unit balls differs substantially at points. However, some basic properties of the extreme points of the balls $B^S_\FM$ and $B^S_\BL$ are worth exhibiting now, as these partially direct the lines of argumentation later.

We start with an observation (cf. \cite{Rao-Roy:1970}, Proposition 2.2, for complex-valued functions on $[0,1]$).
\begin{lemma}\label{lem:trivial extreme points}
	Let $\bullet=\FM,\BL$.
	If $f\in B^S_\bullet$ is such that $|f|=\ind$, then $f\in\ext(B^S_\bullet)$. In particular, if $\bullet=\BL$, then $f=\pm\ind$.
\end{lemma}
\begin{proof}
	Define $S^\pm := \{x\in S: f(x)=\pm1 \}$.
	Let $g\in\BL(S)$ such that both $f+g$ and $f-g$ are in $B^S_\bullet$. Suppose $g\neq 0$ (i.e. $f$ were not extreme). Let $x\in S$ such that $g(x)\neq 0$. If $g(x)>0$ and $x\in S^+$, then $f(x)+g(x)>1$; if $x\in S^-$, then $f(x)-g(x)<-1$. Similar reasoning applies to the case $g(x)<0$, resulting in the conclusion that either $\|f + g\|_\infty>1$ or $\|f-g\|_\infty>1$, contradicting that both $f+g$ and $f-g$ are in $B^S_\bullet$. So, $f$ is extreme.\\
	If $\bullet=\BL$ and $|f|=\ind$, then $\|f\|_\infty=1$ and $|f|_L=0$. So, $f=\pm\ind$.
\end{proof}

A continuous function $f$ on $S$ that satisfies $|f|=\ind$ must be constant on each connected component of $S$, taking the value either $+1$ or $-1$. That is, such $f$ is locally constant. These constitute what we will call {\it trivial extreme points} of $B^S_\bullet$. We shall now be concerned with the {\it non-trivial extreme points}:
\begin{equation}
	\ext_*(B^S_\bullet) := \ext(B^S_\bullet) \setminus \bigl\{ f\in B^S_\bullet: |f|=\ind\bigr\}, \qquad (\bullet=\FM,\BL).
\end{equation}
Note that $\ext_*(B^S_\bullet)=\emptyset$ if $S$ is a trivial space consisting of a singleton. We shall assume that the space $S$ contains at least two elements.

For the $\FM$-norm the conclusion of Corollary \ref{clry:extremes on sphere} (i.e., for an extreme point $f$ of $B^S_\FM$ one has $\max(\|f\|_\infty,|f|_L)=1$) can be strengthened:
\begin{lemma}\label{lem:strengthened norm prop extreme points}
	If $f\in\ext(B^S_\FM)$, then $\|f\|_\infty=1$. If $f\in\ext_*(B^S_\FM)$, then also $|f|_L=1$. 
\end{lemma}
\begin{proof}
	In view of Corollary \ref{clry:extremes on sphere}, $\|f\|_\infty\leq 1$. Suppose, that $\|f\|_\infty<1$. Let $\delta:=1-\|f\|_\infty>0$. Then
	\[
	\|f\pm\delta\ind\|_\infty\leq 1\qquad\mbox{and}\qquad |f\pm\delta\ind|_L \leq |f|_L + \delta|\ind|_L = |f|_L \leq 1.
	\]
	Thus, both $f+\delta\ind$ and $f-\delta\ind\in B^S_\FM$. According to Lemma \ref{lem:char extreme point}, this contradicts that $f$ is an extreme point. Hence, $\|f\|_\infty=1$.
	
	According to the first part, either $|f(x)|=1$ for all $x$, or there exists $x_0\in S$ such that $|f(x_0)|<1$. In the latter case, suppose $|f|_L<1$. Pick $L$ such that $0<L\leq 1-|f|_L$. Let $\delta>0$ such that $|f(x_0)|+\delta\leq1$. Take $h>0$ sufficiently small, such that $(2+\frac{|f|_L}{L})h<\delta$. Define
	\[
	g(x) := \bigl[ h - Ld(x,x_0)\bigr]^+.
	\]
	Then $g\in\BL(S)$, $g\neq 0$ and $|g|_L\leq L$. Moreover, $g(x)=0$ when $d(x,x_0)\geq h/L$ and for $x\in S$ such that $d(x,x_0)<h/L$ one has
	\begin{align*}
		|f(x)\pm g(x)| &\leq  |f(x)-f(x_0)|\ +\ |f(x_0)\pm (h - Ld(x,x_0))|\\
		& \leq |f|_L d(x,x_0) + |f(x_0)|+ h + L d(x,x_0)\\
		& < |f|_L\frac{h}{L} + |f(x_0)| + 2h < |f(x_0)| + \delta \leq 1 .
	\end{align*}
	Thus,
	\[
	\|f\pm g\|_\infty\leq 1\qquad \mbox{and}\qquad |f\pm g|_L \leq |f|_L+|g|_L \leq |f|_L + L \leq 1.
	\]
	So, both $f+g$ and $f-g$ are in $B^S_\FM$, contradicting that $f$ is an extreme point according to Lemma \ref{lem:char extreme point}.
\end{proof}

\begin{remark}
In \cite{Roy:1968}, Theorem 3.1, it has been shown that extreme functions $f$ in the unit ball of $\Lip([0,1])$ for the FM-norm must have $|f'(x)|=1$ almost everywhere with respect to Lebesgue measure, whereas for extreme functions $f$ in the unit ball of $\Lip([0,1]\times[0,1])$ for the FM-norm this no longer holds: the set on which $|f'(x)|=1$ can have arbitrary (but positive) Lebesque measure (less than 1; cf. \cite{Rolewicz:1986}, Theorem 1). This does not violate Lemma \ref{lem:strengthened norm prop extreme points}, since $|f|_L = \mathrm{ess\,sup}_{x\in S} |f'(x)|$.
\end{remark}

For the BL-norm, the situation is different. If $f$ is an extreme point of $B^S_\BL$, then Corollary \ref{clry:extremes on sphere} immediate yields that $|f|_L = 1-\|f\|_\infty$. Interestingly, the following result holds, for which the main inspiration came from \cite{Rao-Roy:1970}, Lemma 2.7. It will be used in the proof of Proposition \ref{E_BL is subset}, which itself is auxiliary to a main result in this paper, namely the generalisation to the BL-norm setting of Johnson's description of a weak-star dense subset of the set of extreme points for the FM-norm, see Corollary \ref{clry:density Johnson-like set BL}. 

\begin{lemma}\label{lem: ext pt attains +-sup}
	If $f\in\ext_*(B^S_\BL)$, then $\inf_{x\in S} f(x) = -\sup_{x\in S} f(x)$. In particular, if $S$ is compact, then $f$ attains both $\|f\|_\infty$ and $-\|f\|_\infty$.
\end{lemma}
\begin{proof}
	Put $M:=\sup_{x\in S} f(x)$, $m:=\inf_{x\in S} f(x)$ and let $\mu:=|M+m|\geq0$. Since $f$ is non-trivial, i.e. $|f|\neq \ind$,  $|f|_L>0$ and $0<\gamma:=\|f\|_\infty = 1-|f|_L<1$. Consequently, $-1<m\leq M<1$ and $\mu<2$. Consider, for $-1< a< 1$ and $b\in \RR$ the functions
	\[
	g_{a,b} := a f + b\ind \in\BL(S).
	\]
	We shall show that if $\mu>0$, then one can find $(a,b)\neq 0$ such that both $f+g_{a,b}$ and $f-g_{a,b}$ are in $B_\BL^S$, while $g_{a,b}\neq 0$. Then Lemma \ref{lem:char extreme point} yields that $f$ cannot be an extreme point of $B^S_\BL$, contradicting the assumption. Hence one must have $M=-m$.
	
	To that end, assume $\mu>0$ and consider the region of interest
	\begin{equation}
		R_\mu := \bigl\{ (a,b)\in\RR^2\,:\, -1<a<1,\ |b| < \smfrac{1}{2}\mu\min(1+a,1-a)\bigr\}
	\end{equation}
	for the parameters $a$ and $b$. $R_\mu$ is defined in such a way that $-|M+m|\leq \frac{2b}{1\pm a}\leq |M+m|$ whenever $(a,b)\in R_\mu$.
	It is an open set containing $(0,0)$. In fact, it is the interior of the closed polyhedral set that is the convex hull of the four points $\pm(1,0)$, $\pm(0,\smfrac{1}{2}\mu)$.
	
	Write $g$ instead of $g_{a,b}$ for notational convenience. Then
	\begin{gather*}
		(1+a)m+ b \leq f+ g \leq (1+a) M + b,\\
		(1-a)m - b \leq f - g \leq (1-a) M - b.
	\end{gather*}
	Hence, given the particular form of $g$,
	\begin{gather*}
		\|f+g\|_\infty = \max\bigl[ (1+a)M+b, -(1+a)m-b \bigr],\\
		\|f-g\|_\infty = \max\bigl[ (1-a)M-b, -(1-a)m+b \bigr],\\
		|f\pm g|_L = (1\pm a)|f|_L = (1\pm a)(1-\|f\|_\infty)
		= (1\pm a)\bigl( 1-\max(M,-m) \bigr).
	\end{gather*}
	Then
	\begin{align}
		\|f+g\|_\BL & = 1 + a + (1+a)\left( \max\left[M+ \frac{b}{1+a}, -m - \frac{b}{1+a}\right] - \max\bigl(M, -m\bigr)  \right)\nonumber\\
		&  = 1+ a + (1+a)\left(\max\left[\frac{2b}{1+a}, -m-M \right] - \frac{b}{1+a} + M - \max\bigl(M, -m\bigr) \right)\nonumber\\
		& = 1 + a - b + (1+a)\left(\max\bigl[ \frac{2b}{1+a},-(M+m)\bigr] -\max\bigl(-(M+m),0\bigr) \right).
	\end{align}
	Similarly, one computes that
	\begin{align}
		\|f-g\|_\BL &= (1-a) + (1-a)\left( \max\bigl[ \frac{2b}{1-a}, M+m\bigr] - \frac{b}{1-a} -m  - \max(M,-m)\right) \nonumber\\
		&= 1-a-b + (1-a)\left(\max\bigl[\frac{2b}{1-a},M+m\bigr] -\max\bigl(M+m,0\bigr) \right).
	\end{align}
	$R_\mu$ has been chosen such that if $(a,b)\in R_\mu$, then
	\begin{equation}
		\|f-g\|_\BL = \begin{cases} 1-a-b, & \mbox{if}\ M+m>0,\\
			1-a+b, & \mbox{if}\ M+m<0 \end{cases}
	\end{equation}
	and
	\begin{equation}
		\|f+g\|_\BL = \begin{cases} 1+a+b, & \mbox{if}\ M+m>0,\\
			1+a-b, & \mbox{if}\ M+m<0 \end{cases}.
	\end{equation}
	Thus, if $M+m>0$ and we require that both $\|f+g\|_\BL\leq 1$ and $\|f-g\|_\BL\leq 1$, then $a=-b$. If $M+m<0$, the the same conditions yield $a=b$. Since $R_\mu$ is open and contains $(0,0)$, in either case $R_\mu$ contains an $(a,-a)\neq 0$ or $(a,a)\neq 0$.\\
	The statement for $S$ compact immediately follows from the general case, since the continuous function $f$ attains its supremum and infimum on the compact set $S$.
\end{proof}

\section{Extension of general Lipschitz functions and extreme points}

The following (elementary) extension result allows limiting our attention to complete spaces $(S,d)$. Let $(\hat{S},\hat{d})$ be the completion of $(S,d)$ and view $S$ as a dense subset of $\hat{S}$.
\begin{prop}\label{prop:extension to completion}
	Every $f\in\Lip(S,d)$ has a unique extension $\hat{f}\in\Lip(\hat{S},\hat{d})$, such that $|\hat{f}|_L=|f|_L$. If $f$ is bounded, then $\hat{f}$ is bounded and $\|\hat{f}\|_\infty = \|f\|_\infty$.
\end{prop}
The proof is left to the reader. According to this result, $\BL(S,d)$ is linearly isometrically isomorphic to $\BL(\hat{S},\hat{d})$. Their respective unit balls and sets of extreme points can hence be identified.

Central in our description of subsets of extreme points $B^S_\BL$ (and  $B^S_\FM$) will be the extension to $S$ of a bounded Lipschitz function on a (finite) subset $P$ -- with the restricted metric. Extension operators for real and vector valued Lipschitz functions have been extensively studied (cf. overviews provided in e.g. \cite{Cobzas-Miculescu-Nicolae:2019,Weaver:1999}, among others). Crucial is what we call the `{\it Metric Tietze's Extension Theorem}', which is stronger than the well-known result by McShane \cite{McShane:1934}, Theorem 1, that states that any Lipschitz function $f$ on $P$ has an extension $F$ to $S$ that has Lipschitz constant that is not greater than that of $f$.
\begin{thrm}[{\bf Metric Tietze's Extension Theorem}]
	\label{thrm:metric Tietze extension}
	Let $(S,d)$ be a metric space and $P\subset S$ a non-empty subset, equipped with the restriction of $d$ to $P$. If $f\in\BL(P)$, then there exists $F\in\BL(S)$ such that $F|_P=f$, $\|F\|_\infty = \|f\|_\infty$ and $|F|_L = |f|_L$.
\end{thrm}
For a full proof, see Appendix \ref{App:proof metric Tietze extension thrm}. One should be aware that extension with preservation of Lipschitz constant is delicate for vector-valued Lipschitz functions. Theorem \ref{thrm:metric Tietze extension} then does not hold for a general codomain (e.g. see \cite{Cobzas-Miculescu-Nicolae:2019}, Section 4.2): there exists spaces $S$ and subsets $P$, such that a Lipschitz function $f:P\to \RR^2$ cannot be extended to $S$ with preservation of the Lipschitz constant. However, in specific cases it may, see e.g. the Kirszbraun-Valentine Theorem, (cf. \cite{Cobzas-Miculescu-Nicolae:2019}, Theorem 4.2.3, p.221), but then by a different extension method. The McShane-type of Lipschitz extension is given by the operator
\begin{equation}\label{eq:def F0}
	\CE_P^{S,0} f(x) := \sup_{p\in P} \,\bigl[ f(p) - |f|_L\,d(p,x)\,\bigr], \qquad x\in S.
\end{equation}
We take the convention that for $P$ a singleton, $|f|_L=0$ for any $f\in\BL(P)$, for the Lipschitz constant is then the supremum in $\RR_+$ of the empty set. Thus, extension from singletons $P=\{x_0\}$ leads to constant functions on $S$.

The extension operator involved in Theorem \ref{thrm:metric Tietze extension} has particularly nice properties with regard to extreme points as we shall show below. It is defined by
\begin{equation}\label{eq:def full extension}
	\CE_P^S:\BL(P)\to\BL(S) : f\mapsto F := \max\bigl(\CE_P^{S,0}( f), -\|f\|_\infty\bigr).
\end{equation}

We start with exhibiting basic properties of the McShane extension operator $\CE_P^{S,0}$.
\begin{prop}\label{prop:props McShane extension}
	Let $P$ be a non-empty subset of $S$, with the metric induced by $S$, and let $f,g\in\BL(P)$, then:
	\begin{enumerate}
		\item[({\it i})] $\CE_P^{S,0}(f)\leq \|f\|_\infty\ind$.
		\item[({\it ii})] $\CE_P^{S,0}(c\ind) = c\ind$ for every $c\in\RR$.
		\item[({\it iii})] If $f\leq g$ and $|f|_L\geq |g|_L$, then
			$\CE_P^{S,0}(f)\leq \CE_P^{S,0}(g)$.
		\item[({\it iv})] If $|f\vee g|_L\geq |f|_L$ and $|f\vee g|_L\geq |g|_L$, then
		$\CE_P^{S,0}(f\vee g) \leq \CE_P^{S,0}(f) \vee \CE_P^{S,0}(g)$.
	\end{enumerate}
\end{prop}
\begin{proof}
	Statements ({\it i}), ({\it ii}) and ({\it iii}) can be derived easily from the definition of the extension operator. For ({\it iv}), let $x\in S$. Then
	\begin{align*}
		\CE_P^{S,0}(& f\vee g)(x)\\
		& =\max\left( \sup_{p\in P:f(p)\geq g(p)} \bigl[ f(p)-|f\vee g|_Ld(x,p) \bigr],\sup_{p\in P:f(p)< g(p)} \bigl[ g(p)-|f\vee g|_Ld(x,p) \bigr]\right)\\
		&\leq \max\left( \sup_{p\in P:f(p)\geq g(p)} \bigl[ f(p)-|f|_Ld(x,p) \bigr],\sup_{p\in P:f(p)< g(p)} \bigl[ g(p)-|g|_Ld(x,p) \bigr]\right)\\
		& \leq \max\bigl(\, \CE_P^{S,0}(f)(x),\, \CE_P^{S,0}(g)(x)\, \bigr).
	\end{align*}
\end{proof}
Note that if the condition in Proposition \ref{prop:props McShane extension} ({\it iv}) holds then equality of Lipschitz constants must hold in at least one of the two inequalities, because in general $|f\vee g|_L\leq \max(|f|_L,|g|_L)$. This observation yields the following
\begin{clry}
	Let $P$ be a non-empty subset of $S$, with the metric induced by $S$.
	If $f\in\BL(P)$ and $g\in\BL(S)$ are such that $\bigl|f\vee g|_P\bigr|_L = |f|_L$ and $\bigl|f\vee g|_P\bigr|_L \geq |g|_L$, then
	\begin{equation}
		\CE_P^{S,0}\bigl(f\vee g|_P\bigr) \vee g \ =\ \CE_P^{S,0}(f) \vee g.
	\end{equation}
\end{clry}
\begin{proof}
	Since $\CE_P^{S,0}(g|_P)\leq g$, we derive from Proposition \ref{prop:props McShane extension} ({\it iv}) that
	\begin{equation}\label{eq:one ineq extension}
		\CE_P^{S,0}\bigl(f\vee g|_P\bigr)\vee g \leq \CE_P^{S,0}(f) \vee g.
	\end{equation}
	Because $\bigl|f\vee g|_P\bigr|_L \leq |f|_L$ and $f\vee g|_P\geq f$, obviously, Proposition \ref{prop:props McShane extension} ({\it iii}) yields
	\[
		\CE_P^{S,0}\bigl(f\vee g|_P\bigr) \geq \CE_P^{S,0}(f).
	\]
	Taking the maximum with $g$ on both sides in the last inequality gives the inequality opposite to \eqref{eq:one ineq extension}.
\end{proof}
Note the following special case of this corollary that is relevant to the extension operator $\CE_P^S$:
\begin{clry}\label{clry:truncation in extension}
	Let $P$ be a non-empty subset of $S$, with the metric induced by $S$. If $f\in\BL(P)$ and $c\in\RR$ are such that $|f\vee c\ind |_L=|f|_L$, then
	\[
		\CE_P^{S,0}\bigl(f\vee c\ind)\vee c\ind = \CE_P^{S,0}(f) \vee c\ind.
	\]
\end{clry}

\begin{remark}\label{rem:extension to completion and McShane}
	If $\overline{P}$ is the closure of $P$ in $S$ (not necessarily complete), then $\CE_P^{\overline{P}}(f)$ is the restriction to $\overline{P}$ of the extension $\hat{f}$ to the completion $\hat{P}$ from Proposition \ref{prop:extension to completion}.
\end{remark}

The extension operators behave as expected under composition:

\begin{prop}\label{prop:repeated extension operators}
	Let $P'\subset S$ and $P\subset P'$ be non-empty subsets, each equipped with the metric induced by $S$. Then $\CE_{P'}^{S,0}\circ \CE_P^{P',0} = \CE_P^{S,0}$ and also $\CE_{P'}^{S}\circ \CE_P^{P'} = \CE_P^{S}$.
\end{prop}
\begin{proof}
	First consider the extension without truncation, defined by $\CE_P^{S,0}$.
	Let $f\in\BL(P)$ and $x\in S$. Then, using that $\bigl|\CE_P^{P',0}(f)\bigr|_L = |f|_L$ according to Theorem \ref{thrm:metric Tietze extension} and the triangle inequality, one obtains:
	\begin{align}
		\bigl(\CE_{P'}^{S,0}\circ \CE_P^{P',0}\bigr)(f) (x) & = \sup_{p'\in P'} \left[ \left( \sup_{p\in P} \bigl[ f(p)- |f|_L d(p',p)\bigr]    \right) - |f|_L d(x,p') \right]\nonumber\\
		& = \sup_{p'\in P'}  \sup_{p\in P}  \left[ f(p) - |f|_L\bigl(d(p',p)+d(x,p')\bigr) \right]\label{eq:intermediate}\\
		&\leq \sup_{p'\in P}  \sup_{p\in P} \bigl[ f(p) - |f|_L d(x,p)\bigr] \ = \ \CE_P^{S,0}(f)(x)\nonumber.
	\end{align}
	On the other hand, equality \eqref{eq:intermediate} can be reformulated, and using $P\subset P'$ one arrives at
	\begin{align*}
		\bigl(\CE_{P'}^{S,0}\circ \CE_P^{P',0}\bigr)(f) (x) & = \sup_{(p,p')\in P\times P'} \left[ f(p) - |f|_L\bigl(d(p',p)+d(x,p')\bigr) \right]\\
		& \geq \sup_{(p,p'):\ p'=p} \left[ f(p) - |f|_L\bigl(d(p',p)+d(x,p')\bigr) \right]\ =\ \CE_P^{S,0}(f)(x).
	\end{align*}

	According to the Tietze Extension Theorem, $\|\CE_P^{P'}(f)\|_\infty = \|f\|_\infty$. If $P$ is not a singleton, because $|\CE_P^{P'}(f)|_L=|f|_L$, there exists a sequence $\bigl((x_n,y_n)\bigr)_{n\in\NN}\in P\times P$, with $x_n\neq y_n$ and $|f(x_n)-f(y_n)|/d(x_n,y_n)\to |f|_L$ as $n\to\infty$. Since $\CE_{P}^{P',0}(f)|_P = f$ and $f\geq -\|f\|_\infty\ind$ on $P$, one obtains by McShane's Theorem that
	\begin{equation}\label{eq:equating Lipschitz constants}
		|f|_L \leq \bigl| \CE_{P}^{P',0}(f) \vee \bigl(-\|f\|_\infty\ind\bigr) \bigr|_L \leq \bigl|\CE_{P}^{P',0}(f) \bigr|_L \leq |f|_L.
	\end{equation}
	Thus, equality of Lipschitz constants must hold everywhere in \eqref{eq:equating Lipschitz constants}.
	If $P=\{x_0\}$ is a singleton, then $\CE_P^{P',0}(f)=f(x_0)\ind$ and equality holds in \eqref{eq:equating Lipschitz constants} trivially. 
	Application of Corollary \ref{clry:truncation in extension} and the first part of the proof now yields
	\begin{align*}
		\bigl(\CE^S_{P'}\circ\CE_{P}^{P'}\bigr)(f) & = \left[ \CE_{P'}^{S,0}\bigl( \CE_{P}^{P'}(f)\bigr)\right]\ \vee\ \bigl[ -\bigl\| \CE_{P}^{P'}(f)\bigr\|_\infty \ind\bigr] \\
		& = \left[ \CE_{P'}^{S,0}\left( \CE_{P}^{P',0}(f) \vee \bigl(-\|f\|_\infty\ind\bigr) \right)\right]\ \vee\ \bigl[ -\bigl\| f\bigr\|_\infty \ind\bigr]\\
		& = \left[ \CE_{P'}^{S,0}\left( \CE_{P}^{P',0}(f)  \right)\right]\ \vee\ \bigl[ -\bigl\| f\bigr\|_\infty \ind\bigr] \\
		& = \ \ \CE_P^{S,0}(f)\ \vee \ \bigl[ -\bigl\| f\bigr\|_\infty \ind\bigr]\ =\ \CE_P^S(f).
	\end{align*}
\end{proof}

Farmer \cite{Farmer:1994} considered extreme points of the Lipschitz dual $S^\#$, i.e. the vector space $\Lip_0(S)$ of Lipschitz functions that vanish at a fixed distinguished point $e\in S$, equipped with the $|\cdot|_L$ seminorm, that becomes a norm on this subspace of $\Lip(S)$. He showed that the McShane extension operator $\CE_P^{S,0}$ maps extreme points of the unit ball in $P^\#$ to extreme points of the unit ball in $S^\#$, for closed subsets $P$ of $S$ (cf. \cite{Farmer:1994}, Lemma 2, p. 810). Hidden in the techniques of proof in \cite{Johnson:1975} for the case of compact and connected $S$, a similar result could be observed for $\BL(S)$ with the FM-norm (stated precisely in Appendix \ref{App:Norming sets FM norm}). The following result establishes this result for the BL-norm. Its novel method of proof is applicable to the FM-norm as well with small modification (see Appendix \ref{App:Norming sets FM norm}) without restricting conditions of connectedness nor compactness. Hence, through the identification of $(\BL(S),\|\cdot\|_\FM$ with $(\Lip_0(S),|\cdot|_L)$ detailed in Appendix \ref{app:isomormism FM-case}, it provides a novel proof for that result as well.

\begin{thrm}\label{thrm:extension non-trivial extreme points}
	Let $P$ be a subset of $S$ with at least two points, equipped with the restriction of $d$ as metric. Then $\mathcal{E}_P^S$ maps $\ext_*(B_{\BL}^P)$ into $\ext_*(B_{\BL}^S)$.
\end{thrm}
\begin{proof}
	Let $\overline{P}$ be the closure of $P$ in $S$. According to Remark \ref{rem:extension to completion and McShane}, $\CE_P^{\overline{P}}$ {\it identifies} $B^P_\BL$ with $B^{\overline{P}}_\BL$ and $\ext_*(B^P_\BL)$ with $\ext_*(B^{\overline P}_\BL)$. In view of Proposition \ref{prop:repeated extension operators} it then suffice to show that $\CE_{\overline{P}}^S$ maps $\ext_*(B^{\overline P}_\BL)$ into  $\ext_*(B^S_\BL)$. That is, without loss of generality we can assume that $P$ is closed.
	
	Let $f\in \ext_*(B_{\BL}^P)$ and $F:=\mathcal{E}_P^S(f)$. Theorem \ref{thrm:metric Tietze extension} yields $F|_P=f$, $\|F\|_\infty=\|f\|_\infty$ and $|F|_L=|f|_L$. Since $f\neq\pm\ind$, $F\neq\pm\ind$. So, it suffices to show that $F$ is an extreme point of $B_{\BL}^S$. Suppose that $G\in \BL(S)$ is such that both $F + G$ and $F-G$ are in $B_{\BL}^S$. We have to show that $G=0$.
	
	We have $f\pm G|_P=(F\pm G)|_P\in B_{\BL}^P$ and $f\in\ext(B_{\BL}^P)$, hence by Lemma \ref{lem:char extreme point}: $G|_P=0$ .
	Define $M_F^-:=\{x\in S:F(x)=-\|F\|_\infty\}$ and $\tilde{H}(s,p):=\frac{H(p)-H(s)}{d(s,p)}$ for $H\in \BL(S)$, $s\in S\setminus P,$ $p\in P$. Note that $|\tilde{H}(s,p)|\leq|H|_L$.
	
	Let $x\in(M_F^-\cup P)^c$. Then $F(x)>-\|F\|_\infty=-\|f\|_\infty$, hence, by definition of $\mathcal{E}_P^S$: $F(x)=\sup_{p\in P}[f(p)-|f|_Ld(x,p)]$. Pick $(p_n)\subset P$ with
	\begin{align}
		F(x)=\lim_{n\to\infty}[f(p_n)-|f|_Ld(x,p_n)].\label{hoi1}
	\end{align}
	Note that the sequence $(f(p_n))$ is bounded. Then, the sequence $(d(x,p_n))$ is bounded by \eqref{hoi1}. Thus, there exists a subsequence $(p_{n_k})$ of $(p_n)$ such that $(d(x,p_{n_k}))$ converges, and we have $\lim_k d(x,p_{n_k})\geq d(x,P)>0$, since $P$ could be assumed closed. Together with \eqref{hoi1}, this implies (recall $F|_P=f$):
	\begin{align}
		\lim_{k\to\infty}\tilde{F}(x,p_{n_k})=\lim_{k\to\infty}\frac{F(p_{n_k})-F(x)}{d(x,p_{n_k})}=|f|_L. \label{hoi2}
	\end{align}
	We have ($F|_P=f$, $G|_P=0$)
	\[
	\|F\pm G\|_\infty\geq \|(F\pm G)|_P\|_\infty=\|f\|_\infty
	\]
	Also, $f\in\ext(B_{\BL}^P)$, so $\|f\|_{\BL}=1$ (otherwise $\|f\|_{\BL}<1$ so $f\pm \eps\in B_{\BL}^P$ for some $\eps>0$).
	So for all $p\in P$, we have
	\[
	|f|_L=1-\|f\|_\infty\geq 1-\|F\pm G\|_\infty\geq |F\pm G|_L\geq |\hat{F}(x,p)\pm\hat{G}(x,p)|,
	\]
	thus $|\hat{F}(x,p)+|\hat{G}(x,p)||\leq |f|_L$.
	In particular, (recall \eqref{hoi2} and $G|_P=0$)
	\[
	\frac{|G(x)|}{d(x,p_{n_k})}=|\hat{G}(x,p_{n_k})|\leq |f|_L-\hat{F}(x,p_{n_k})\quad \longrightarrow\ 0\quad \text{as } k\to\infty.
	\]
	Since $(d(x,p_{n_k}))$ is bounded, we conclude that $G(x)=0$.
	
	It remains to show that we have $G|_{M_F^-\cap P^c}=0$. Let $x\in M_F^-$, i.e., $F(x)=-\|F\|_\infty$. We have
	\[
	|F\pm G|_L\geq |(F\pm G)|_P|_L=|f|_L,
	\]
	so (recall $\|f\|_{\BL}=1$):
	\begin{align}
		\|f\|_\infty =1-|f|_L\geq 1-|F\pm G|_L\geq \|F\pm G\|_\infty\geq |F(x)\pm G(x)|.\label{hoi3}
	\end{align}
	Now, if $G(x)>0$, then $|F(x)-G(x)|=\|F\|_\infty+G(x)>\|f\|_\infty$, while if $G(x)<0$, we have $|F(x)+G(x)|=\|F\|_\infty+G(x)>\|f\|_\infty$, in both cases contradicting \eqref{hoi3}. Thus $G(x)=0$.
\end{proof}

\section{Norming and dense sets of extreme points}

Our interest in the extreme points of the unit ball of Banach spaces of Lipschitz functions originates from the definition of useful norms on $\CM(S)$ by means of embedding measures in the dual of such Banach spaces and Proposition \ref{prop:norm attained in extreme point}. The latter reduces the task of computing $\|\mu\|_\bullet^*$, $\bullet=\FM$ or $\BL$, to determining $\sup\bigl\{\langle\mu,f\rangle\colon f\in\ext(B^S_\bullet)\bigr\}$, provided $\mu$ is separable. 

The full set of extreme points of $B_\bullet^S$ cannot be conveniently described for this purpose. The K$^2$-M$^3$-R Theorem, Theorem \ref{thrm:K2-M3-R}, indicates what type of sets to look for in order to take the supremum over a smaller set of extreme points. Recall the natural pairing $\langle\cdot,\cdot\rangle$ between $\CM(S)$ and $\BL(S)$, defined by integration, as in \eqref{eq:pairing functions-measures}. Proposition \ref{prop:dual Ms is BL}, which asserts that $(\BL(S),\|\cdot\|_\BL)$ is the dual space of $\overline{\Mol}(S)_\BL=\CM_s(S)$, yields that the  $\sigma(\BL(S),\overline{\Mol}(S)_\BL)$-weak topology defined by the above pairing is the weak$^*$-topology on $\BL(S)$. According to the Banach-Alaoglu Theorem (\cite{Conway:1990}, Theorem V.3.1), $B^S_\bullet$ is compact in this locally convex Hausdorff vector space topology. Thus, according to Theorem \ref{thrm:K2-M3-R} ({\it iv}), a subset $A$ of $\ext(B_\bullet^*)$ is norming for $\|\cdot\|_\bullet^*$ if and only if any of the conditions $({\it i})$-$({\it iii})$ in the K$^2$-M$^3$-R Theorem holds.

We provide sets $E_\bullet^S$ of extreme points of which the weak$^*$-closed convex hull is $B^S_\bullet$, thus proving property ({\it i}). Lemma \ref{lem:weakstar dense subsets unit ball} provides the technique to do so. As a consequence one gets the weak$^*$-density in the extreme points (which is property ({\it ii})). As `by-product' of this approach, it is possible to reduce attention to the dense subset $\Mol(S)$ of $\CM(S)$ of measures with {\it finite support}. This will allow to define `small' norming sets $E^S_\bullet$, defined by means of the Lipschitz extension operators $\CE^S_P$. Johnson, in the setting of connected and compact metric spaces in \cite{Johnson:1975}, aimed at providing a set $J^S_\FM$ of extreme points in $B^S_\FM$ that was `as large as possible', possibly having a different objective in mind. 

For $f\in \BL(S)$, define 
\[
M_f := \{x\in S\colon |f(x)| = \|f\|_\infty \},
\]
which is non-empty in the case $S$ is compact. Johnson defined
\begin{align}
J^S_\FM := \bigl\{ f\in B^S_\FM\colon& \|f\|_\infty = 1,\ \exists P_f\subset S\ \mbox{finite},\ P_f\neq\emptyset:\label{def: J-set FM}\\
& \bigr(\; x\in S\setminus M_f\ \Rightarrow\ \exists p\in P_f: |f(x)-f(p)|=d(x,p)\;\bigr)\  \bigr\}\nonumber
\end{align}
(cf. the set $\mathcal{A}$ defined in \cite{Johnson:1975}, Proposition 1.1). Below we shall now define `our' small sets $E^S_\bullet$, but also provide for the case of the $\BL$-norm a novel `Johnson-like' large set $J^S_\BL$ of extreme points of $B^S_\BL$. Motivation for the definition of $J^S_\BL$ comes from Lemma \ref{lem: ext pt attains +-sup}.

Thus, we introduce the subset of $B^S_\bullet$:
\begin{equation}
	\hat{E}^S_\bullet := \bigcup_{P\subset S,\ \mathrm{finite}} \CE_P^S\bigl(\ext_*(B^P_\bullet)\bigr),\qquad \bullet=\FM,\BL.
\end{equation}
According to Theorem \ref{thrm:extension non-trivial extreme points} (for BL-norm) and Theorem \ref{prop:extension non-trivial extr points} (for FM-norm), $\hat{E}^S_\bullet$ consists of non-trivial extreme points of $B^S_\bullet$, obtained by McShane-type Lipschitz extension. Note that $\ext_*(B^P_\bullet)=\emptyset$ when $P$ is a singleton. So, the union runs over all finite sets $P$ with at least two elements. Expand the set $\hat{E}^S_\BL$ by adding the trivial extreme points: 
\begin{equation}
	E^S_\BL := \hat{E}^S_\BL \ \cup\ \bigl\{ \pm\ind \},
\end{equation}
For the FM-norm, there are non-trivial extreme functions that cannot be reached by the Lipschitz extension operators $\CE_P^S$, see Appendix \ref{App:Norming sets FM norm}, but which must be included for the norming property to hold. These need to be added separately, apart from the trivial extreme points. So, in this case: 
\begin{equation}\label{def:E-FM}
	E^S_\FM := \hat{E}^S_\FM \ \cup\ \bigl\{f\in B^S_\FM\colon |f|=\ind \} \ \cup\ \bigl\{ h_P\colon P\subset S\ \mbox{finite},\ P\neq\emptyset \bigr\},
\end{equation}
where
\begin{equation}\label{def:hP main text}
	h_P(x) := (-\ind)\,\vee\, \sup_{p\in P} \bigl[ 1 - d(x,p) \bigr],\qquad x\in S.
\end{equation}
Lemma \ref{lem:non-extended extreme points} shows that the functions $h_P$ are extreme points of $B^S_\FM$.
Note that these functions $h_P$ cannot be reached by extension from a Lipschitz function on $P$, because $(h_P)|_P=\ind$ and $\CE_P^S(\ind) = \ind$ (cf. Proposition \ref{prop:props McShane extension} ({\it ii})).

The Johnson-like set $J^S_\BL$ is defined similar to \eqref{def: J-set FM}:
\begin{align}
	\hat{J}^S_\BL &:=  \bigl\{ f\in B^S_\BL\colon
	\|f\|_\BL=1,\  f(M_f)=\{\|f\|_\infty, -\|f\|_\infty \},\label{def:hat J-set BL} \\
	&\qquad\quad \exists P_f\subset S\ \mbox{finite},\ P_f\neq\emptyset:\nonumber\\
	&\qquad\quad \bigr(\; x\in S\setminus M_f\ \Rightarrow\ \exists p\in P_f: |f(x)-f(p)|=(1-\|f\|_\infty)d(x,p)\;\bigr)\  \bigr\}\nonumber\\
	J^S_\BL &:= \hat{J}^S_\BL\ \cup\ \bigl\{\pm\ind\}\label{def: J-set BL}
\end{align}
Note that in the case of the FM-norm, the trivial extreme points, satisfying $|f|=\ind$, are contained in $J^S_\FM$ as formulated in \eqref{def: J-set FM}, because $M_f=S$ for such $f$. Hence, the condition with $P_f$ is trivially satisfied. For the BL-norm, the trivial extreme points $f=\pm\ind$ also have $M_f=S$, but $f(M_f)$ does not consist of both two values $\pm\|f\|_\infty$. It is necessary though, to have the trivial extreme points `on board'.
Note too, that the functions $h_P$, which were added to $\hat{E}^S_\FM$ to get $E^S_\FM$, are included in $J^S_\FM$.



\begin{prop}\label{prop:E in J}
	$E^S_\bullet\subset J^S_\bullet$ for $\bullet=\BL$ or $\FM$.
\end{prop}
\begin{proof}
	{\it Case of FM-norm.}\ If $f\in E^S_\bullet$ is trivial, $|f|=\ind$ and hence $M_f=S$. Thus, we need to check only the conditions not involving $P_f$. One has $\|f\|_\infty =1$, so $f\in J^S_\FM$.\\
	If $f=h_P$ for some $\emptyset\neq P\subset S$ finite, then $\|f\|_\infty=1$ and $M_f= P\cup \{x\in S\colon f(x)=-1\}$. One may take $P_f=P$ and then the second condition in the definition is satisfied.\\
	If $f\in\hat{E}^S_\FM$ is non-trivial, then $\|f\|_\infty$ and $|f|_L=1$, according to Lemma \ref{lem:strengthened norm prop extreme points}. Moreover,  $f=\CE^S_P(f^*)$ for some $f^*\in\ext_*(B^P_\FM)$ and $P\subset S$ finite, with $|P|\geq 2$. By the Tietze Extension Theorem, $|f^*|_L=|f|_L=1$ and it is clear from \eqref{eq:def F0} and \eqref{eq:def full extension} that $f$ satisfies the second condition in the definition, using $P_f=P$.\\
	{\it Case of BL-norm.}\ If $f\in E^S_\BL$ is trivial, $|f|=\ind$ and hence $M_f=S$. Thus, again we need to check only the conditions not involving $P_f$. In this case, $f=\pm\ind$, according to Lemma \ref{lem:trivial extreme points}, which are included in $J^S_\BL$ by construction.\\
	If $f\in \hat{E}_\BL^S$, then $f=\CE_P^S(f^*)$ for some $f^*\in\ext_*(B_\BL^P)$ with $P$ finite. From \eqref{eq:def F0} and \eqref{eq:def full extension} it is immediate that $\CE_P^S(f^*)$ satisfies the last condition appearing in the definition of $\hat{J}_\BL^S$ with $P_f:=P$. The other two conditions follow from Corollary \ref{clry:extremes on sphere} and Lemma \ref{lem: ext pt attains +-sup} applied to $f^*\in\ext_*(B_\BL^P)$. Hence, $f=\CE_P^S(f^*)\in \hat{J}_\BL^S\subset J_\BL^S$.
\end{proof}

The next result will clarify the relationship between $J^S_\BL$ and $\ext(B^S_\BL)$. To that end we need the following lemma, which is a reformulation of \cite{Rao-Roy:1970}, Lemma 2.1 for real numbers. 

\begin{lemma}\label{lemma2.1royrao}
	Let $x,y\in\RR$ and let $|x|+|y|=1$. If $\alpha,\beta\in\RR$ are such that $|x\pm\alpha|+|y\pm\beta|\leq 1$, then either $xy=\alpha=\beta=0$, or
	\[
	xy\neq 0,\quad |\alpha|\leq\min\{|x|,|y|\}\quad\text{ and }\quad \alpha\frac{|x|}{x}+\beta\frac{|y|}{y}=0.
	\]
\end{lemma}
\begin{proof}
	See \cite{Rao-Roy:1970}, Lemma 2.1. 
\end{proof}

We can only prove that $J^S_\BL$ (and $J^S_\FM$, see Appendix \ref{App:Norming sets FM norm}) consist of extreme points if $S$ is connected. The prove for the BL-norm is remarkably more complicated than that for the FM-norm. It is unclear whether $J^S_\bullet$ always contains non-extremal points when $S$ is not connected. Example \ref{ex:point in Johson not extreme} gives an example of a non-extremal point when $S$ is finite.

\begin{prop}\label{E_BL is subset}
	If $S$ is connected, then  $J_{\bullet}^S\subset\ext(B_{\bullet}^S)$ for $\bullet=\FM$ or $\BL$.
\end{prop}
\begin{proof}
	Assume that $S$ is connected.\\
	{\it Case of FM-norm:}\ The argument is essentially that in \cite{Johnson:1975}, proof of Proposition 1.1, with removal of the compactness condition. Suppose that $f\in J^S_\FM$ and let $P_f$ be the corresponding set $P$. Suppose that $g\in \BL(S)$ is such that both $f+g$ and $f-g$ are in $B^S_\FM$. We shall show that $g=0$ and conclude by means of Lemma \ref{eq:cond extreme point} that $f$ is extreme.\\
	If $x\in M_f$, then $f(x)\in\{\pm1\}$ and $1\geq\|f\pm g\|_\infty\geq |f(x)\pm g(x)|$. Therefore $g(x)=0$. If $x\in S\setminus (M_f\cup P_f)$, then there exists $p_x\in P_f$ such that $\tilde{f}(x,p_x):= \frac{f(p_x) - f(x)}{d(x,p_x)} = \pm 1$. Define $\tilde{g}(x,p_x)$ similarly. Then for both choices of sign,
	\[
	1\geq \|f\pm g\|_\FM \geq |f\pm g|_L \geq |\tilde{f}(x,x_p)\pm \tilde{g}(x,x_p)| = |1\pm \tilde{g}(x,x_p)|.
	\]
	So, $\tilde{g}(x,x_p)=0$ and $g(x)=g(x_p)$. Therefore,
	\[
	g(S) = g(M_f) \cup g(S\setminus(M_f\cup P_f))\cup g(P_f) = g(P_f).
	\]
	In particular, $g(S)$ is finite. Since $g$ is continuous and $S$ is connected, $g(S)$ consists of a singleton.\\
	If $M_f\neq\emptyset$ (e.g. if $S$ is compact), then $0\in g(M_f)\subset g(S)$. So $g=0$. If $M_f=\emptyset$, then also $\{0\}\subset g(S)$, by contradiction: suppose that $0\not\in g(S)$. Pick a sequence $(x_n)\subset S$ such that $|f(x_n)|\to\|f\|_\infty = 1$. Then
	\[
	|f|+|g| = \max\bigl( |f+g|, |f-g|\bigr) \leq \max\bigl( \|f+g\|_\infty, \|f-g\|_\infty\bigr) \leq 1
	\]
	Consequently,
	\[
	0\leq|g(x_n)|\leq 1- |f(x_n)| \to 0\qquad \mbox{as}\ n\to\infty.
	\]
	So, $|g(x_n)|\to 0$. However, $g(S)$ is finite. Thus, $g(x_n)=0$ must hold for sufficiently large $n$, which implies that $0\in g(S)$. We arrived at a contradiction. Hence, $g=0$ and we conclude that $f$ is extreme.
	
	{\it Case of BL-norm:}\ Suppose that $f\in J_{\BL}^S$ and $f\pm g\in B_{\BL}^S$ for some $g\in \BL(S)$. We show that $g=0$, i.e.  $f\in\ext(B_{\BL}^S)$.	
	Recall that $J_{\BL}^S=\hat{J}_{\BL}^S\cup\{\pm\ind\}$.
	If $f=\pm\ind$, it is immediate that $g$ must be zero.
	In the other case, we have $f\in \hat{J}_{\BL}^S$.
	For $s,t\in S$, $s\neq t$, define $\tilde{f}(s,t):=\frac{f(s)-f(t)}{d(s,t)}$ and  $\tilde{g}(s,t):=\frac{g(s)-g(t)}{d(s,t)}$.
	Let $P_f$ be a finite subset as in the definition of $\hat{J}_{\BL}^S$, and let $s\in M_f^c\cap P_f^c$. Set
	$\gamma:=\|f\|_\infty$ and note that $\gamma\in(0,1)$. Now, let $p_s\in P_f$ be such that $|\tilde{f}(s,p_s)|=1-\gamma$. Such $p_s$ exists according to the definition of $J_{\BL}^S $.
	Then, for all $x\in M_f$ we have
	\begin{align*}
		\begin{cases}
			&|f(x)|+|\tilde{f}(s,p_s)|=\gamma+1-\gamma=1,\\
			&|f(x)\pm g(x)|+|\tilde{f}(s,p_s)\pm \tilde{g}(s,p_s)|\leq\|f\pm g\|_\infty+|f\pm g|_L=\|f\pm g\|_{\BL}\leq 1.
		\end{cases}
	\end{align*}
	Note that $f(x)\neq0\neq \tilde{f}(s,p_s)$ since $\gamma\in(0,1)$, so Lemma \ref{lemma2.1royrao} yields
	\begin{align}
		g(x)\frac{|f(x)|}{f(x)}+\tilde{g}(s,p_s)\frac{|\tilde{f}(s,p_s)|}{\tilde{f}(s,p_s)}=0, \qquad \text{ for all } x\in M_f. \label{heyy}
	\end{align}
	Keeping $s$ (and $p_s$) fixed and varying $x$, this implies
	\begin{align}
		g(x)=-\frac{\tilde{g}(s,p_s)|\tilde{f}(s,p_s)|}{\gamma\tilde{f}(s,p_s)}f(x)=cf(x), \qquad \text{ for all } x\in M_f,\label{heyyy}
	\end{align}
	where $c:=-\frac{\tilde{g}(s,p_s)|\tilde{f}(s,p_s)|}{\gamma\tilde{f}(s,p_s)}=\frac{g(x)}{f(x)}$ is constant.
	Now, \eqref{heyy} becomes
	\[
	c\gamma+\tilde{g}(s,p_s)\frac{|\tilde{f}(s,p_s)|}{\tilde{f}(s,p_s)}=0,
	\]
	so $\tilde{g}(s,p_s)=\frac{-c\gamma}{1-\gamma}\tilde{f}(s,p_s)$ (recall that $|\tilde{f} (s,p_s)|=1-\gamma$). Therefore,
	\begin{align}
		g(s)-g(p_s)=\frac{-c\gamma}{1-\gamma}(f(s)-f(p_s)). \label{heyhoi}
	\end{align}
	Since $f\in J_{\BL}^S$, there exist $x_{\pm}\in M_f$ with $f(x_{\pm})=\pm\|f\|_\infty=\pm\gamma$.
	For the moment, assume the following:
	\begin{align}
		\begin{cases}&\text{There exist } n\in \NN \text{ and } \{x_1,\ldots,x_n\}\subset S, \text{ such that } x_1=x_{-}, x_n=x_{+}\\
			&\text{and }
			g(x_{i+1})-g(x_i)=\frac{-c\gamma}{1-\gamma}(f(x_{i+1})-f(x_i)) \text{ for }
			i=1,\ldots,n-1.\end{cases} \label{aanname}
	\end{align}
	We will prove \eqref{aanname} below, but we finish the proof first, using this assumption.
	We have
	\begin{align*}
		2c\gamma &=cf(x_n)-cf(x_1)
		\overset{\eqref{heyyy}}{=}g(x_{n})-g(x_{1})\\
		&=\sum_{i=1}^{n-1}\big(g(x_{i+1})-g(x_i)\big)
		=\sum_{i=1}^{n-1}\frac{-c\gamma}{1-\gamma}\big(f(x_{i+1})-f(x_i)\big)\\
		& =\frac{-c\gamma}{1-\gamma}(f(x_n)-f(x_1))\\
		& =\frac{-c\gamma}{1-\gamma}2\gamma.
	\end{align*}
	So either $c=0$ or $\gamma=0$ or  $1=\frac{-\gamma}{1-\gamma}$. The last two options cannot hold since $\gamma\in(0,1)$.
	We conclude that $c=0$.
	Now, \eqref{heyyy} implies that $g=cf=0$ on $M_f\supset\{x_+,x_-\}$ and from \eqref{heyhoi}, we derive that $g(s)\in g(P_f)$ for all $s\in S\setminus P_f$.
	Thus, $\{0\}\subset g(S)\subset\{0\}\cup g(P_f)$, so $g(S)$ is finite and contains zero.
	Since $S$ is connected and $g$ is continuous, $g(S)$ is connected. Also, it is  finite and $0\in g(S)$, so $g(S)=\{0\}$, i.e. $g=0$.  We conclude that $f\in\ext(B_{\BL}^S)$.
	
	It remains to prove \eqref{aanname}.
	To this end, define for $x\in S$:
	\begin{align}
		M_x:=\{s\in S:&\exists n\in \NN,  \exists x_1,\ldots,x_n\in S \text{ such that }  x_1=x, x_n=s,\notag \\
		&g(x_{i+1})-g(x_i)=\frac{-c\gamma}{1-\gamma}(f(x_{i+1})-f(x_i))
		\text{ for } 1\leq i\leq n-1
		\}.\label{def Mx}
	\end{align}
	We write $x\sim y$ if and only if $y\in M_x$.
	Clearly, if $x_{+}\in M_{x_{-}}$, we have proved \eqref{aanname}.
	
	We first show that $\sim$ is an equivalence relation.
	By taking $n=1$ in the definition of $M_x$ we see that trivially $x\in M_x$, proving reflexivity.
	For symmetry, suppose that $y\in M_x$, with corresponding $\{x_1,\ldots,x_n\}$ as in the definition of $M_x$. Define $\bar{x}_i:=x_{n-i+1}$ for $i=1,\ldots,n$. Then $\bar{x}_1=x_n=y$, $\bar{x}_n=x_1=x$ and
	\begin{align*}
		g(\bar{x}_{i+1})-g(\bar{x}_i)=-\Big(g(x_{n-i+1})-g(x_{n-i}) \Big)&=-\Big(\frac{-c\gamma}{1-\gamma}(f(x_{n-i+1})-f(x_{n-i})) \Big)\\
		&=\frac{-c\gamma}{1-\gamma}(f(\bar{x}_{i+1})-f(\bar{x}_i)),
	\end{align*}
	so  $x\in M_y$.
	For transitivity, suppose that $x\in M_y$ with corresponding $\{x_1,\ldots,x_n\}$ and $y\in M_z$ with corresponding $y_1,\ldots,y_m$. Define $N:=n+m$ and $\bar{x}_i:=x_i$ for $1\leq i\leq n$, $\bar{x}_i:=y_{i-n}$ for $n+1\leq i\leq N$. Note that $\bar{x}_n=\bar{x}_{n+1}=y$ and note that all the properties are satisfied to conclude that $x\in M_z$.
	We conclude that $\sim$ is an equivalence relation. Consequently, any two equivalence classes $M_x, M_y$ are  either equal or disjoint.
	
	We show that each $M_x$ is a closed subset of $S$. Suppose that $(s_n)\subset M_x$ and  $s_n\to s\in S$. By continuity of $f$, we have $f(s_n)\to f(s)$.
	Since $f\in J_{\BL}^S$, we have for each $n\in \NN$ either $s_n\in M_f$, so $|f(s_n)|=\gamma$, or $|f(s_n)-f(p)|=(1-\gamma)d(s,p)$ for some $p\in P_f$. Write $P_f=\{p_1,\ldots,p_k\}$ with $p_i\in S$ and $k\in\NN$. Define
	\begin{align*}
		&I_i:=\{n\in \NN: |f(s_n)-f(p_i)|=(1-\gamma)d(s_n,p_i)\},\\
		&I^{\gamma}:=\{n\in\NN: f(s_n)=\gamma\},\quad
		I^{-\gamma}:=\{n\in\NN: f(s_n)=-\gamma\}.
	\end{align*}
	It holds that $\NN=\cup_{i=1}^k I_i\cup I^{\gamma}\cup I^{-\gamma}$, so at least one of the sets on the right-hand side must be infinite. Let $I$ be such an infinite set on the right-hand side and write $I=\{n_1,n_2,\ldots\}$ with $n_j<n_{j+1}$ for all $j\in\NN$, so that $(s_{n_j})$ is a subsequence of $(s_n)$.
	
	If $I=I_i$  for some $i\in\{1,\ldots,k\}$, then
	\[
	|f(s)-f(p_i)|=\lim_{j\to\infty}|f(s_{n_j})-f(p_i)|=\lim_{j\to\infty}(1-\gamma)d(s_{n_j},p_i)=(1-\gamma)d(s,p_i),
	\]
	using that $f(s_{n_j})\to f(s)$ and  $d(s_{n_j},s)\to 0$. By \eqref{heyhoi}, it now  follows that $s\in M_{p_i}$ (take $n=2$, $x_1=p_i$, $x_2=s$ in the definition \eqref{def Mx}, note that $s\in M_{p_i}$ trivially holds when $s=p_i$).
	Moreover, $s_{n_j}\in M_{p_i}$ by \eqref{heyhoi}, so $ p_i\sim s$ and $p_i\sim s_{n_j}$ for all $j\in\NN$, implying $s_{n_j}\sim s$. Also, $s_{n_j}\in M_x$ for all $j\in\NN$, so $x\sim s_{n_j}$.  Consequently, $x\sim s$, i.e. $s\in M_x$.
	
	The other cases are $I=I^{\pm\gamma}$. In these cases, we have $f(s_{n_j})=\pm \gamma$ for all $j\in\NN$, so $f(s)=\lim_{j\to\infty}f(s_{n_j})=\pm \gamma=f(s_{n_j})$. Therefore, $s,s_{n_j}\in M_f$, and \eqref{heyyy} yields $g(s_{n_j})=g(s)=\pm c\gamma$. Now, trivially we have $g(s_{n_j})-g(s)=\frac{-c\gamma}{1-\gamma}(f(s_{n_j})-f(s))$ as both sides are zero, and it follows that $s_{n_j}\in M_s$. Thus, $s\sim s_{n_j}$ and $x\sim s_{n_j}$ for all $j\in\NN$, implying $x\sim s$, i.e. $s\in M_x$.
	We conclude that $M_x$ is closed.
	
	Next, observe that
	\begin{align}
		S=\cup_{i=1}^k M_{p_i}\cup M_{x_+}\cup M_{x_-}.\label{heyyyy}
	\end{align}
	For any $s\in S\setminus M_f$, \eqref{heyhoi} holds for some $p_s\in\{p_1,\ldots,p_k\}=P_f$, hence $s\in M_{p_s}$. For $s\in M_f$ we have $f(s)=\pm\gamma=f(x_{\pm})$ and $g(s)=\pm c\gamma=g(x_{\pm})$, so taking $n=2,x_1=x_{\pm}, x_2=s$ in the definition of $M_{x_{\pm}}$, we see that $s\in M_{x_{\pm}}$. So indeed, \eqref{heyyyy} holds.
	
	Our goal was to show that $x_{+}\in M_{x_{-}}$. Suppose that the latter is not true. Then $M_{x_+}\cap M_{x_-}=\emptyset$. Let $J:=\{i\in 1,\ldots,k\}:p_i\notin M_{x_+}\}$ and let $M:= \cup_{i\in J} M_{p_i}\cup M_{x_-}$.
	Then $M$ is a finite union of closed sets, hence closed. Moreover, $M\cap M_{x_+}=\emptyset$ and $x_-\in M_{x_-}\subset M$, so $M\neq \emptyset$. Also, $M_{x_+}$ is closed,  $x_+\in M_{x_+}$ and $S=M\dot\cup M_{x_+}$,
	contradicting the connectedness of $S$.
	We conclude that $x_+\in M_{x_-}$, justifying \eqref{aanname}.
\end{proof}

The following result is essential in establishing the weak$^*$-density of $E^S_\bullet$ in $\ext(B^S_\bullet)$, ($\bullet=\FM$ or $\BL$).


\begin{prop}\label{prop:E is norm-determining}
	Let $(S,d)$ be a metric space and $\bullet=\FM$ or $\BL$. For every $\mu\in\Mol(S)$, there exists $f^\mu_\bullet\in E_\bullet^S$ such that $\|\mu\|_\bullet^*=\pair{\mu}{f^\mu_\bullet}$.
\end{prop}
\begin{proof}
	Let $\mu\in\Mol(S)$ and put $P:=\supp(\mu)=\{s_1,\dots,s_n\}$ with all $s_i$ distinct. So, $\mu=\sum_{j=1}^n a_j\delta_{s_j}$, with $0\neq a_j\in\RR$. Let $P$ inherit the metric $d$ from $S$. View $\mu$ as an element of $\BL(S)^*$, or as an element $\mu|_P$ of $\BL(P)^*$, given by the same expression.
	
	For any $f\in\BL(S)$ the restriction $f|_P\in\BL(P)$ and $\|f|_P\|_\infty\leq \|f\|_\infty$, $|f|_P|_L\leq |f|_L$. So, $\|f|_P\|_\bullet\leq \|f\|_\bullet$ and consequently, $\{f|_P:P\in B^S_\bullet\}\subset B^P_\bullet$. Moreover, $\pair{\mu|_P}{f|_P} = \pair{\mu}{f}$. Hence,
	\begin{equation}\label{eq:comparison bullet-norms P and S}
	 \|\mu\|_\bullet^* = \sup\bigl\{ \pair{\mu|_P}{f|_P}: f\in B^S_\bullet\bigr\}\ \leq\ \|\mu|_P\|_\bullet^*
	\end{equation}
 	Since $(P,d)$ is a finite metric space, it is Polish and Proposition \ref{prop:norm attained in extreme point} yields the existence of $f^*\in\ext(B^P_\bullet)$ such that $\|\mu|_P\|_\bullet^* = \pair{\mu|_P}{f^*}$.
	
	If $f^*$ is non-trivial, then $F^*:=\CE_P^S(f^*)\in E_\bullet^S$ is such that $\pair{\mu}{F^*} = \pair{\mu|_P}{F^*|_P}= \pair{\mu|_P}{f^*} = \|\mu|_P\|^*_\bullet$. Since $\|\mu\|_\bullet^*\geq \pair{\mu}{F^*}$, equation \eqref{eq:comparison bullet-norms P and S} yields that $\pair{\mu}{F^*}=\|\mu\|_\bullet^*$, so we may take $f^\mu_\bullet := F^*$.
	
	If $f^*$ is trivial, i.e. $|f^*|=1$, then in case $\bullet=\BL$, $f^*=\pm\ind$ and one can take $f^\mu_\BL:=\pm\ind\in E_\BL^S$. The case $\bullet=\FM$ is slightly more involved. Consider $P^+:=\{p\in P: f^*(p)=1\}$ and $h:=h_{P^+}\in\ext(B^S_\FM)$ as defined in Lemma \ref{lem:non-extended extreme points} as function on $S$. Since for any $p\in P\setminus P^+$, $f^*(p)=-1$ and $|f^*|_L\leq 1$, one has for any $p\in P\setminus P^+$ and $p'\in P^+$
	\[
	\frac{|f^*(p) - f^*(p')|}{d(p,p')} = \frac{2}{d(p,p')} \leq 1.
	\]
	So $d(p, P^+)\geq 2$. Lemma \ref{lem:non-extended extreme points}, applied to extension from $P^+$ to $P$, shows that $h|_P = f^*$. $h$ is included in $E^S_\FM$ by construction.
\end{proof}

We can now prove the main result and its immediate consequences.

\begin{thrm}\label{thrm:main result dense set extreme points}
	Let $(S,d)$ be a metric space and $\bullet=\FM$ or $\BL$. Then $E_\bullet^S$ is a weak$^*$-dense subset of $\ext(B^S_\bullet)$. If, additionally, $S$ is compact, then $E_\bullet^S$ is $\|\cdot\|_\infty$-dense in $\ext(B^S_\bullet)$.
\end{thrm}

\begin{proof}
	We shall employ Lemma \ref{lem:weakstar dense subsets unit ball} and  Theorem \ref{thrm:K2-M3-R} to obtain the desired result. First we apply Lemma \ref{lem:weakstar dense subsets unit ball} to the situation where $X=\overline{\Mol}(S)_\bullet$ and $F=\Mol(S)$, which is dense in $X$ for the $\|\cdot\|_\bullet^*$-norm. Proposition \ref{prop:dual Ms is BL} yield that $X^*=\BL(S)$. Proposition \ref{prop:E is norm-determining} yields that $E_\bullet^S$ is determining the $\|\cdot\|_\bullet^*$-norm on $F$. Lemma \ref{lem:weakstar dense subsets unit ball} allows to conclude that $B^S_\bullet$ is the weak$^*$-closed convex hull of $E^S_\bullet$.\\
	Moreover, by the Banach-Alaoglu Theorem (cf. \cite{Conway:1990}), $B^S_\bullet$ is compact for the weak$^*$-topology. Apply Theorem \ref{thrm:K2-M3-R} to the locally convex space $\BL(S)$, equipped with the weak$^*$-topology. Statement ({\it i}) has been proven above. So, ({\it ii}) holds. That is, $E^S_\bullet$ is weak$^*$-dense in $\ext(B^S_\bullet)$.
	
	Assume now that $S$ is compact. Let $f\in\ext(B^S_\bullet)$. By the first part, there exists a net $(f_\alpha)_{\alpha\in A}$ in $E^S_\bullet$ such that $f_\alpha\to f$ weak$^*$. Thus, $f_\alpha(x)\to f(x)$ for every $x\in S$, by applying $\delta_x\in \Mol(S)$. Since $S$ is metric and compact, it is totally bounded. Let $\eps>0$. Then there exist $x_1,\dots x_n$ such that $S=\bigcup_{i=1}^n B(x_i,\eps/3)$. For each $i$ one can find $\alpha_i\in A$ such that for all $\alpha \succeq \alpha_i$, $|f_\alpha(x_i)-f(x)|<\eps/3$. Then there also exist $\alpha_0\in A$ such that the same inequality holds for all $1\leq i\leq n$ and all $a\succeq\alpha_0$. Pick $x\in S$. There exists $x_i$ such that $x\in B(x_i,\eps/3)$. Then for all $\alpha\succeq\alpha_0$,
	\begin{align*}
		|f\alpha(x)-f(x)|\ &\leq\ |f_\alpha(x) - f_\alpha(x_i)|\ +\  |f_\alpha(x_i) - f(x_i)|\ + \ |f(x_i)-f(x)|\\
		& \leq\ |f_\alpha|_L d(x,x_i)\ +\ \smfrac{1}{3}\eps\ +\ |f|_L d(x,x_i)\ <\ \eps.
	\end{align*}
	We conclude that $\|f_\alpha - f\|_\infty\to 0$.
\end{proof}

In view of the arguments at the start of this section, we immediately obtain from the K$^2$-M$^3$-R Theorem that the sets $E^S_\bullet$ are norming on $\CM_s(S)$:
\begin{clry}\label{clry:E is norming}
	Let $\mu\in\CM_s(S)$. Then $\|\mu\|_\bullet^* = \sup\{\langle \mu,f\rangle\colon f\in E^S_\bullet\}$, for $\bullet=\FM,\BL$.
\end{clry}

Moreover, we obtain weak$^*$-density of the large set $J^S_\BL$ in the extreme points, provided $S$ is connected.

\begin{clry}\label{clry:density Johnson-like set BL}
	If $(S,d)$ is a connected metric space, then the set $J^S_\bullet$ consists of extreme points $J^S_\bullet$ of $B^S_\bullet$ and is  weak*-dense in $\ext(B^S_\bullet)$ ($\bullet=\FM,\BL$).
\end{clry}
\begin{proof}
	Combine Proposition \ref{E_BL is subset} and Theorem \ref{thrm:main result dense set extreme points}.
\end{proof}

\section{Considerations for the case of a finite space}

In this section we take a closer look at the interpretation of the results of the previous sections for the case when $S$ is a metric space with finitely many points. This seems a highly specific setting, but it is relevant for multiple reasons. Firstly, for a general metric space the relevant norming and dense set $E^S_\bullet$ is obtained mostly by applying the Lipschitz extension operator to (non-trivial) extreme points of the unit ball of functions on a finite metric space $P$. It is then a natural question to see to what extent the latter are determined by the extreme points of a ball for a space $P'\subset P$. Secondly, highly relevant in view of the question of computing FM- and BL-norms of measures, the norm of a molecular measure $\mu$ with support $\supp(\mu)=P$, which is finite, is determined by the integrals against $\mu$ of the functions in the finite set $\ext(B_\bullet^P)$ (see \cite{Hille-Theewis:2022}, Section 6). Below we limit attention to relatively novel case of the BL-norm. Although $S$ is finite, we shall keep talking about functions on $S$.

Theorem \ref{thrm:main result dense set extreme points} provides density of
\begin{equation}\label{def disconnected E_BL}
E^S_\BL =\bigcup_{P\subset S,\ \mathrm{finite}} \CE_P^S\bigl(\ext_*(B^P_\BL)\bigr)\cup\ \bigl\{ f\in B^S_\BL: |f|=1\bigr\}
\end{equation}
in $\ext(B^S_\BL)$. When $S$ is finite, we can take $P=S$ on the right-hand side and \eqref{def disconnected E_BL} becomes trivial (knowing Theorem \ref{thrm:extension non-trivial extreme points}).
Still, let us study \eqref{def disconnected E_BL} a bit further. One may wonder, whether it is possible to reach any non-trivial extremal function by applying the extension operator to a suitable non-trivial extremal function defined on a {\it proper} subset $P$ of $S$.

One can quickly see that this is not the case. Since $B^S_\bullet$ is symmetric around $0$, also $\ext(B^S_\bullet)$ is symmetric around $0$ (cf. Corollary \ref{clry:symmetric sets and extreme points}). If $f=\CE_P^S(F)$ for some $F\in\ext_*(B^P_\bullet)$, then $-f$ is not of this form, generally. This motivates to study yet another extension operator:
\[
\bar{\CE}_P^S(f):=-\CE_P^S(-f)=\min\left(\inf_{p\in P}[f(p)+|f|_Ld(p,\cdot)],\|f\|_\infty\right).
\]
It is easily seen that $\bar{\CE}_P^S$ enjoys many of the same properties as $\CE_P^S$. In particular, $\bar{\CE}_P^S$ maps $\BL(P)$ into $\BL(S)$ and preserves supremum norms and Lipschitz constants.

The properties of $\CE_P^S$ expressed in Proposition \ref{prop:repeated extension operators} and Theorem \ref{thrm:extension non-trivial extreme points} indicate that for any $P\subset P' \subset S$:
\begin{align}\label{eq:inclusion inductive}
	\CE_{P'}^S \circ\CE_{P}^{P'}\bigl(\ext_*(B_\bullet^{P})\bigr)=\CE_{P}^{S}\bigl(\ext_*(B_\bullet^{P})\subset \ext_*(B_\bullet^{S}).
\end{align}
Since $\ext_*(B_\bullet^S)=-\ext^*(B_\bullet^S)$, \eqref{eq:inclusion inductive} also holds with $\CE_P^S$ replaced by $\bar{\CE}_P^S$.
Thus for any finite $S$ containing at least two points we have
\begin{align}\label{ext pts induction}
\bigcup\limits_{x\in S}\mathcal{E}_{S\setminus\{x\}}^S\big(\ext_*(B_{\BL}^{S\setminus\{x\}})\big)\cup  \bar{\CE}_{S\setminus\{x\}}^S\big(\ext_*(B_{\BL}^{S\setminus\{x\}})\big)
\subset \ext_*(B_\BL^{S}).
\end{align}
A natural question is then whether equality holds in \eqref{ext pts induction}.

For the moment, assume that indeed we have equality in \eqref{ext pts induction}
and suppose that $n:=|S|\in\NN_{\geq 3}$.
By induction, we would then have 
\begin{align}\label{wrong hypothesis}
\ext_*(B_\BL^{S})=\bigcup\limits_{\substack{x,y\in S,\\ x\neq y}} \bigcup\{ T_n\circ\ldots\circ T_3(\ext_*(B_{\BL}^{\{x,y\}})): \,& T_i=\CE_{P_{i-1}}^{P_i} \text{ or }\, T_i=\bar{\CE}_{P_{i-1}}^{P_i}, \notag\\[-.5cm]
&\{x,y\}\subset P_i\subset P_{i+1}\subset S \text{ and } |P_i|=i \}.
\end{align}
Now, we could construct all extreme points explicitly, since we `only' have to take all combinations of extensions of the extreme points for $P=\{x,y\}$, which are explicitly given by
\begin{equation}\label{ext pts two points}
\ext_*(B_{\BL}^{\{x,y\}})=\{(x,y)\mapsto \pm (\frac{d(x,y)}{d(x,y)+2},\frac{-d(x,y)}{d(x,y)+2})\}.
\end{equation}
This can be see easily by identifying $\BL(P)\cong \RR^2$ and drawing $B^P_\BL$. However, the next example shows that in general the other inclusion in \eqref{ext pts induction} does NOT hold and \eqref{wrong hypothesis} fails.  Thus, one cannot obtain all extreme points via the inductive procedure described. 
In addition, the example shows that `$P\subset S$' could not have been replaced by `$P\subsetneq S$'  in \eqref{def disconnected E_BL}.

\begin{example}\label{example no induction}
Let $S$ be finite with $|S|\geq 3$.
If equality in \eqref{ext pts induction}, hence \eqref{wrong hypothesis}, were true, then for all $f\in\ext_*(B_\BL^S)$ there would exist
$x,y \in S$ with $f(x)=\|f\|_\infty, f(y) = -\|f\|_\infty$ and $2\|f\|_\infty=d(x,y)(1-\|f\|_\infty)$. Indeed, the property holds for $S=P$, $|P|=2$ by \eqref{ext pts two points} and obviously, it is preserved by $\CE_P^S$ and $\bar{\CE}_P^S$.

Now let $S:=\{0, 1.5, 2.5, 4\} \subset \RR$ and identify $f\in\BL(S)$ with $(f(0),f(1.5),f(2.5),f(4))\in\RR^4$. Using the expression for $\ext(B_{\BL}^S)$ in terms of linear constraints in \cite{Hille-Theewis:2022}, one can (numerically) verify that $f:=(0.5, -0.25, 0.25, -0.5)$ is an extreme point of $B_\BL^S$. However, $2\|f\|_\infty=1\neq 2= d(0,4)(1-0.5)$, so it does not satisfy the aforementioned property.

We conclude that equality in \eqref{ext pts induction} does not hold and \eqref{wrong hypothesis} is false in general.
\end{example}

Now we turn our attention to the `Johnson-like' set $J^S_\BL$. For finite $S$, the property that defines $J^S_\BL$ holds for all extreme points, not only for those in $E^S_\BL$ (cf. Proposition \ref{prop:E in J}):  
\begin{lemma}\label{lem:johnson property finite S}
Let $S$ be finite, $f\in \ext(B_{\BL}^S)$ and $x\in S\setminus M_f$. Then there exists
$y\in S\setminus\{x\}$ such that $|f(x)-f(y)|=|f|_Ld(x,y)$.
\end{lemma}
\begin{proof}
Suppose the claim is false. Then we have $|f(x)-f(y)|<|f|_Ld(x,y)$ for all $y\in S\setminus\{x\}$. 
Also, $|f(x)|<\|f\|_\infty$, so
\[
\alpha:=(\|f\|_\infty-|f(x)|)\wedge \min_{y\in S\setminus\{x\}}[|f|_Ld(x,y)-|f(x)-f(y)|]>0.
\]
Put $g:=\ind_{\{x\}}\alpha\in \BL(S)$ (since $S$ is finite). We show that $\|f\pm
g\|_{\BL}\leq 1$, contradicting the fact that $f\in \ext(B_{\BL}^S)$.
We have
\begin{align*}
\|f\pm g\|_\infty =\max\{\|f\|_\infty,|f(x)\pm \alpha|\}&\leq
\max\{\|f\|_\infty,|f(x)|+\alpha\}\\
&\leq\max\{\|f\|_\infty,|f(x)|+\|f\|_\infty-|f(x)|\}\\
&=\|f\|_\infty.
\end{align*}
Furthermore,
\begin{align*}
|f\pm g|_L&=|f|_L\vee\max_{y\in S\setminus \{x\}}\frac{|(f\pm g)(x)-(f\pm
g)(y)|}{d(x,y)}\\
&=|f|_L\vee\max_{y\in S\setminus \{x\}}\frac{|f(x)-f(y)\pm \alpha|}{d(x,y)}\\
&\leq |f|_L\vee\max_{y\in S\setminus \{x\}}\frac{|f(x)-f(y)|+|f|_Ld(x,y)-|f(x)-f(y)|}{d(x,y)}\\
&=|f|_L.
\end{align*}
Thus $\|f\pm g\|_{\BL}\leq \|f\|_\infty+|f|_L=1$, finishing the proof.
\end{proof}

If $S$ is a finite metric space then it is disconnected, so the proof of $J^S_\BL\subset \ext(B_\BL^S)$ does not apply. Also, the definition of $J_\BL^S$ now boils down to $J^S_\BL=\{f\in B^S_B\colon \|f\|_\BL=1,\ f(M_f)=\{\|f\|_\infty,-\|f\|_\infty\}\,\}$, since one can take $P_f=S$ for the finite set $P_f$ in the definition. However, Lemma \ref{lem:johnson property finite S} makes one wonder whether there is a variation of the defining property of $J^S_\BL$ that does characterize $\ext(B_\BL^S)$. E.g. is the set
\begin{align}\label{def finite johnson set}
\{f\in\BL(S): \|f\|_\BL=1, f(M_f)=&\{\pm\|f\|_\infty\}, \forall s\in M_f^c\; \exists s\in S\setminus\{s\}\notag\\
&\text{ s.t. } |f(p)-f(s)|=(1-\|f\|_\infty)d(p,s)\}
\end{align}
contained in $\ext(B_\BL^S)$? Then we would have equality by Lemma \ref{lem:johnson property finite S}, Lemma \ref{lem: ext pt attains +-sup} and Corollary \ref{clry:extremes on sphere}, i.e. a full characterization of $\ext(B_\BL^S)$.
Unfortunately, the next example proves the contrary.

\begin{example}\label{ex:point in Johson not extreme}
Let $S:=\{0, 1.5, 2, 4\}$ and identify $f\in\BL(S)$ with $(f(0),f(1.5),f(2),f(4))\in\RR^4$. Then $f:=(0.5, -0.25, 0, -0.5)$ is an element of the set \eqref{def finite johnson set}. However, one can (numerically) verify that it is not an extreme point of $B_\BL^S$ (e.g. as in Example \ref{example no induction} using \cite{Hille-Theewis:2022}).

We conclude that for finite $S$, the properties from Lemma \ref{lem:johnson property finite S}, Lemma \ref{lem: ext pt attains +-sup} and Corollary \ref{clry:extremes on sphere} are necessary conditions to be in $\ext(B_\BL^S)$, but not sufficient conditions.
\end{example}


\begin{appendix}

\section{Proof of Metric Tietze Extension Theorem}
\label{App:proof metric Tietze extension thrm}

We could not find a reference for the Metric Tietze Extension Theorem in the literature that provides a full proof. See \cite{Weaver:1999}, Theorem 2.5.6 though, which does not provide full details. The result is fundamental for our results. Therefore, we include a full detailed proof here. 

\begin{proof} {\it (Metric Tietze Extension Theorem, Theorem \ref{thrm:metric Tietze extension}).}\\
	Let $F_0=\mathcal{E}_{P}^{S,0}(f)$ and $F=\mathcal{E}_P^S(f)$ be as defined in \eqref{eq:def F0} and \eqref{eq:def full extension}.
	First, we show that
	$F|_P=F_0|_P=f$. Let $p_0\in P$.
	For all $p\in P$, we have
	\[
	f(p)-|f|_Ld(p,p_0)\leq f(p)-(f(p)-f(p_0))=f(p_0)=f(p_0)-|f|_Ld(p_0,p_0),
	\]
	thus $F_0(p_0)=f(p_0)$. Moreover, $f(p_0)\geq -\|f\|_\infty$ yields  $F(p_0)=F_0(p_0)=f(p_0)$, proving that $F|_P=F_0|_P=f$.
	
	As a consequence, we have $\|F\|_\infty\geq\|f\|_\infty$. For the other inequality, note that $f(p)-|f|_Ld(p,x)\leq \|f\|_\infty$ for any $x\in S$ and $p\in P$, so $F_0\leq \|f\|_\infty$. From the definition of $F$, it follows that $-\|f\|_\infty\leq F\leq \|f\|_\infty$, i.e. $\|F\|_\infty\leq \|f\|_\infty$.
	Hence, $\|F\|_\infty= \|f\|_\infty$.
	
	
	We have $|F|_L\geq |F|_P|_L= |f|_L$, provided that $F$ is Lipschitz continuous. It remains to show  that the latter holds and  $|F|_L\leq |f|_L$.
	Let $x,y\in S$ and $p\in P$. Applying the triangle inequality twice yields
	\[
	f(p)-|f|_Ld(p,y)-|f|_Ld(x,y)\leq f(p)-|f|_Ld(p,x)\leq f(p)-|f|_Ld(p,y)+|f|_Ld(x,y).
	\]
	Taking the supremum over $p\in P$, we obtain
	\[
	F_0(y)-|f|_Ld(x,y)\leq F_0(x)\leq F_0(y)+|f|_Ld(x,y).
	\]
	Therefore, $|F_0(x)-F_0(y)|\leq |f|_Ld(x,y)$ and $F_0$ is Lipschitz continuous.
	It follows that  $F=\max(F_0,-\|f\|_\infty)$ is Lipschitz continuous with $|F|_L\leq |F_0|_L\leq |f|_L$ (\cite{Dudley:1966}, Lemma 4). We conclude that $|F|_L=|f|_L$.
\end{proof}

\section{Isomorphism between $(\BL(S),\|\cdot\|_\FM)$ and the metric dual space}
\label{app:isomormism FM-case}
	
Farmer \cite{Farmer:1994} considered the Lipschitz dual $S^\#$ of the metric space $(S,d)$, which was introduced by Lindenstrauss \cite{Lindenstrauss:1964}. It consists of the vector space $\Lip_0(S)$ of all functions in $\Lip(S)$ that vanish at a fixed chosen `distinguished point' $e\in S$, equipped with $|\cdot|_L$ as norm. We exhibit here a sequence of isomorphisms that identify $(\BL(S,d),\|\cdot\|_{\FM,d})$ with a such a space $(\Lip_0(S^+,d^+),|\cdot|_{L,d^+})$, linearly isometrically, order theoretically and algebraically. These isomorphisms can be found in \cite{Weaver:1999}, Section 1.7.

Let $(S,d)$ be a metric space. Define 
\[
d'(x,y) := d(x,y)\wedge 2.
\]
\begin{lemma}\label{lem:first reduction}
	The identity map is a linear isomertic, lattice and algebraic isomorphism between $(\BL(S,d),\|\cdot\|_{\FM,d})$ and $(\BL(S,d'),\|\cdot\|_{\FM,d'})$.
\end{lemma}
\begin{proof}
	We leave the elementary check to the reader. See e.g. \cite{Weaver:1999}, Proposition 1.7.1.
\end{proof}

For a metric space $(S,d')$ of diameter at most two, i.e. 
\[
\mathrm{diam}(S,d') := \sup\bigl\{ d'(x,y)\colon, x,y\in S\,\}\leq 2,
\]
like the space in Lemma \ref{lem:first reduction}, one can add a remote distinguished point $e$, not in $S$, to obtain $S^+:= S\cup\{e\}$ with metric
\[
d^+(x,y) := \begin{cases} 
	d(x,y), & \mbox{if}\ x,y\in S,\\
	1, & \mbox{if}\ x=e,\ \mbox{or}\ y=e,\ x\neq y,\\
	0, & \mbox{if}\ x=e=y.
	\end{cases}	
\]
$Lip_0(S^+,d^+)$ is then the space of Lipschitz funtions on $S^+$ for $d^+$ that vanish at $e$.
\begin{lemma}
	If $\mathrm{diam}(S,d')\leq 2$, then $(\BL(S),\|\cdot\|_{\FM,d'})$ is linearly isometrically, order theoretically and algebraically isomorphic to $(\Lip_0(S^+,d^+), |\cdot|_{L,d^+}\,)$.
\end{lemma}
\begin{proof}
	The isomorphism is given by extending $f\in \BL(S)$ by 0 at $e$. See \cite{Weaver:1999}, Theorem 1.7.2.
\end{proof}

The two results show that the set $B^S_\FM$ remains the same set of functions in its image in $\Lip_0(S^+,d^+)$, except that each function has been trivially extended to the added distinguished point $e$. This unit ball becomes the unit ball in $\Lip_0(S^+,d^+)$. Its convex structure does not change, so extreme points correspond.

\section{An alternative proof and auxiliary result for the FM-norm}
\label{App:Norming sets FM norm}

Through the sequence of isomorphisms provided above in Appendix \ref{app:isomormism FM-case}, the case of the FM-norm reduces essentially to already established results on the Lipschitz dual space, e.g. those by Farmer \cite{Farmer:1994}. In his proofs he uses the characterisation of extreme points, {\it loc.cit}, Theorem 1. An equivalent of this is not available for the BL-norm. However, the equivalent of Theorem \ref{thrm:extension non-trivial extreme points}, which is central in the construction of the small sets of extreme points $E^S_\bullet$, does hold for the FM-norm, see \cite{Farmer:1994}, Lemma 2. Below we shall present a proof of this result, without using Farmer's characterisation theorem, but instead follow the line of reasoning in the proof of Theorem \ref{thrm:extension non-trivial extreme points}. It requires only minor modifications.

Recall that $\ext_*(B^P_\bullet)=\emptyset$ if $P$ is a singleton.

\begin{thrm}\label{prop:extension non-trivial extr points}
	Let $P$ be a subset of $S$ with at least two elements, equipped with the restriction of $d$ as metric. Then, $\CE_P^S$ maps $\ext_*(B^P_\FM)$ into $\ext_*(B^S_\FM)$.
\end{thrm}
\begin{proof}
	A similar argument as at the start of the proof of Theorem \ref{thrm:extension non-trivial extreme points} shows that without loss of generality we may assume $P$ to be closed.
	
	Theorem \ref{thrm:metric Tietze extension} yields that $\CE_P^S$ maps $B^P_\FM$ isometrically into $B^S_\FM$. Let $f\in \ext_*(B^P_\FM)$ and put $F:=\CE_P^S(f)$. According to Theorem \ref{thrm:metric Tietze extension} and Lemma \ref{lem:strengthened norm prop extreme points} one has $\|F\|_\infty = \|f\|_\infty=1$ and $|F|_L=|f|_L=1$. Since $|f|\neq\ind$ on $P$ and $F|_P=f$, $|F|\neq \ind$ on $S$. Therefore, it remains to show that $F$ is extreme. It suffices to show that if $F+G$ and $F-G$ are both in $B^S_\FM$ for some $G\in\BL(S)$, then $G=0$ (cf. Lemma \ref{lem:char extreme point}).
	
	First, $(F\pm G)|_P=f\pm G|_P\in B^P_\FM$, with $G|_P\in\BL(P)$. $f$ is an extreme point of $B^P_\FM$, so $G|_P=0$. Define $M_F^-:=\{x\in S: F(x)=-1\}$. For $H\in B^S_\FM$, define
	\begin{equation}
		\hat{H}(s,p) := \frac{H(p)-H(s)}{d(s,p)},\qquad s\in S\setminus P,\ p\in P.
	\end{equation}
	Then $|\hat{H}(s,p)|\leq 1$. Let $x\in (M_F^-\cup P)^c$. There exist $p_n\in P$, $n\in\NN$, such that
	\begin{equation}\label{eq:convergence}
		f(p_n) - d(x,p_n)\ \uparrow\ F_0(x) = F(x)\quad \mbox{as}\ n\to\infty,
	\end{equation}
	Because $f$ is bounded and \eqref{eq:convergence} holds, the sequence $(d(x,p_n))_n$ must be bounded. Since $P$ is assumed to be closed and $x\not\in P$, $\sup_n d(x,p_n)\geq \inf_n d(x,p_n)>0$. Let $(p_{n_k})$ be a subsequence such that $(d(x,p_{n_k}))$ converges (necessarily to a non-zero limit). Then, \eqref{eq:convergence} yields
	\[
	\hat{F}(x,p_{n_k})\ \to\ 1\quad\mbox{as}\ k\to\infty.
	\]
	Since both $F+G$ and $F-G$ are in $B^S_\FM$, we get
	\[
	1 \geq \|F\pm G\|_\FM \geq |F\pm G|_L \geq |\hat{F}(x,p) \pm \hat{G}(x,p)|\qquad\mbox{for all}\ p\in P,
	\]
	which implies $\bigl|\hat{F}(x,p) + |\hat{G}(x,p)|\bigr|\leq 1$ for all $p$. Since $G|_P=0$, we find that for sufficiently large $k$,
	\[
	\frac{|G(x)|}{d(x,p_{n_k})}\ =\ |\hat{G}(x,p_{n_k})|\ \leq\ 1 - \hat{F}(x,p_{n_k})\ \to 0\quad\mbox{as}\ k\to\infty.
	\]
	Because the sequence $(d(x,p_{n_k}))_k$ has non-zero limit, $G(x)=0$.
	If $x\in M_f^-$, then $F(x)=-1$ and
	\[
	1\geq \|F \pm G\|_\FM\geq \|F\pm G\|_\infty \geq |F(x) \pm G(x)| = |1\mp G(x)|.
	\]
	We conclude that $G(x)=0$ in this case too. So $G=0$ and $f$ is an extreme point.
\end{proof}

Recall the definition \eqref{def:hP main text} of the functions $h_P$. It was remarked that these cannot be reached by Lipschitz extension operators. The following lemma summarizes fundamental properties of these functions $h_P$.  
\begin{lemma}\label{lem:non-extended extreme points}
	Let $P$ be a closed, proper and non-empty subset of $S$.
	Then $h_P\in\ext(B^S_\FM)$. If there exists $x\in S\setminus P$ such that $d(x,P)<2$, then $h_P\in\ext_*(B^S_\FM)$. Otherwise, $h_P$ is a trivial extreme point.
\end{lemma}
\begin{proof}
	By construction, $h_P\in B^S_\FM$ and $h_P(p)=1$ for all $p\in P$. Let $g\in\BL(S)$ be such that $h_P+g$ and $h_P-g$ are both in $B^S_\FM$. Then $g(p)=0$ for all $p\in P$. Put $M_h^-:=\{x\in S: h_P(x)=-1\}$. Then $g=0$ on $M^-_h$ too, since both $h_P(x)+g(x)\geq -1$ and $h_P(x)-g(x)\geq -1$ for $x\in M^-_h$. If $S\setminus(M^-_h\cup P)=\emptyset$, then $|h_P|=\ind$, so $h_P$ is a trivial extreme point of $B^S_\FM$.
	
	Suppose $S\setminus(M^-_h\cup P)\neq\emptyset$ and pick $x\in (M^-_h\cup P)^c$. Since $P$ is closed, $d(x,P)>0$. Moreover,  there exist $p_n\in P$, $n\in\NN$, such that $h_P(x) = \lim_{n\to\infty} 1- d(x,p_n)$. Since for any $p\in P$ one has $g(p)=0$ and $h_P(p)=1$, and $h_P\pm g\in B^S_\FM$,
	\begin{equation}\label{eq:basic estimate}
		\frac{|h_P(x)\pm g(x) - 1|}{d(x,p)} \leq |h_P \pm g|_L \leq 1\qquad \mbox{for all}\ p\in P.
	\end{equation}
	Let $\eps>0$. There exists $N\in\NN$ such that $|h_P(x)-1+d(x,p_n)|< \eps$ for all $n\geq N$.
	From \eqref{eq:basic estimate} one derives that
	\[
	\bigl| d(x,p_n)\pm g(x)\bigr| \leq \bigl|h_P(x)-1+d(x,p_n)\bigr| + \bigl|h_P(x)\mp g(x)-1\bigr| < \eps + d(x,p)
	\]
	for all $n\geq N$ and $p\in P$. Taking $p=p_n$, we obtain $|g(x)|< \eps$. Since $\eps>0$ was taken arbitrarily, $g(x)=0$.
	Hence $g=0$ on $S$ and $h_P\in\ext(B^S_\FM)$, according to Lemma \ref{lem:char extreme point}.
	
	Finally, assume that there exists $x_0\in S\setminus P$ such that $d(x_0,P)<2$. Then $-1<h_P(x_0)<1$, so $h_P$ cannot be a trivial extreme point. If for every $x\in S\setminus P$, $d(x,p)\geq 2$ for all $p\in P$. Then $h_P(x)=-1$ for $x\in S\setminus P$, while $h_P(p)=1$ for every $p\in P$. Hence, $h_P$ is a trivial extreme point.
\end{proof}

\end{appendix}

\end{document}